\documentclass[11pt,DIV=14]{scrartcl}
\usepackage{amsmath,amssymb,amsfonts}
\usepackage[colorlinks=true,citecolor=red,urlcolor=blue,linkcolor=blue,pdfborder={0 0 0}]{hyperref}
\usepackage{stmaryrd}

\newtheorem{theorem}{Theorem}[section]%
\newtheorem{lemma}[theorem]{Lemma}
\newtheorem{remark}[theorem]{Remark}

\allowdisplaybreaks

\providecommand{\proof}{\textbf{Proof: }}
\providecommand{\close}{\nolinebreak\mbox{$\blacktriangleleft$}}

\newcommand{\R}{\mathbb{R}}

\providecommand{\xx}{\mathbf x }

\providecommand{\N}{\mathbb{N}}

\providecommand{\eps}{\varepsilon}
\providecommand{\BB}{{\mathbf B }}
\providecommand{\DD}{{\mathbf D }}
\providecommand{\EE}{{\mathbf E }}
\providecommand{\HH}{{\mathbf H }}
\providecommand{\JJ}{{\mathbf J }}
\providecommand{\Th}{{\mathcal{T}_h }}

\title%
{An Energy Stable Discontinuous Galerkin Time-Domain Finite Element Method in Optics and Photonics}
\author{%
A.~Anees\thanks{%
Dept.~of Mathematics and Statistics, University of Agriculture,
Faisalabad 38000, Pakistan, asadanees@uaf.edu.pk},
L.~Angermann\thanks{%
Dept.~of Mathematics, Clausthal University of Technology, Erzstr.~1, D-38678 Clausthal-Zellerfeld, Germany,
lutz.angermann@tu-clausthal.de}}

\begin{document}
\maketitle

\begin{abstract}
In this paper, a time-domain discontinuous Galerkin (TDdG) finite element method
for the full system of Maxwell's equations in optics and photonics is investigated,
including a complete proof of a semi-discrete error estimate.
The new capabilities of methods of this type are to efficiently model linear and
nonlinear effects, for example of Kerr nonlinearities.
Energy stable discretizations both at the
semi-discrete and the fully discrete levels are presented. In particular,
the proposed semi-discrete scheme is optimally convergent in the spatial variable
on Cartesian meshes with $Q_k$-type elements,
and the fully discrete scheme is conditionally stable with respect
to a specially defined nonlinear electromagnetic energy.
The approaches presented prove to be robust and allow the modeling of optical problems
and the treatment of complex nonlinearities as well as geometries of various physical systems
coupled with electromagnetic fields.
\end{abstract}



\section{Introduction}
There is great interest in developing time-domain discontinuous Galerkin (TDdG) methods for the
full system of Maxwell's equations in optics and photonics, for instance to design optical devices
with higher complexity. One of the most famous and general problems is
the third order Kerr-type nonlinear model.
A few excellent and efficient schemes are available on developing
the finite difference time-domain (FDTD) methods, to solve the nonlinear Maxwell's equations
in Kerr media (e.g., \cite{Joseph:94, Ziolkowski:93, Bokil:18, Jia:19, Wang:21}).
There are also many studies of TDFEMs for Maxwell's equations for considering the flexibility
of finite element methods for complex domains and materials.
Recent advancements and more references on TDFEMs and TDdG for Maxwell's equations
with Kerr-type nonlinearity can be found in some recent reviews such as
\cite{Fisher:07, Huang:21, Angermann:20a, Peng:20,
Abraham:19b, Jiang:19, Abraham:19a, Abraham:19c}.

In the past few years, TDdG have gotten considerable attention and are
being employed to a wide range of problems in optics and photonics.
To the authors' knowledge, few mathematical proofs for the convergence
of the discontinuous Galerkin method when applied to Maxwell's equations
were given in the papers \cite{Hesthaven:02, Jiang:19}.
These methods allowed a comparatively easy handling of elements of various types and shapes,
irregular non-conforming meshes and even locally varying polynomial degrees.
There are still accessibly few analysis, error estimates, and simulation results
by employing TDdG available for the system of Maxwell's equations with Kerr-type nonlinearity.
Moreover, the results were given mostly for the 1D case using TDdG schemes.

Our prime object is to develop an energy stable numerical scheme in this paper
that can preserve the stability relation at the semi-discrete and fully discrete levels.
Energy preserving schemes are robust since
they are able to maintain and preserve the shape and phase of the waves accurately
after long-term numerical simulations.
Moreover, error estimates are also presented at the semi-discrete and fully discrete levels.
In this paper, we extend our results about semi-discrete conforming mixed finite element methods
\cite{Angermann:16b}, \cite{Angermann:18a}
and fully discrete conforming finite element methods
\cite{Angermann:18b}, \cite{Angermann:19g}, \cite{Angermann:20a}
to a discontinuous Galerkin discretization of the Maxwell's equations with nonlinearities.
For the sake of simplicity, we will present the result in 2D,
and analogous result are obtained for 3D.
At the end of this brief and by no means complete overview on the related literature
it should be mentioned that, in the course of preparing this work,
the paper \cite{Lyu:21} was published with the same intention and comparable results,
the authors of which probably were unaware of the first author's PhD thesis \cite{Anees:20}.
An essential part of the presented paper is an revised outcome of \cite{Anees:20},
in which the ideas and results about the proposed method were formulated for the first time.

In the spatial discretization we use $Q_k$-type elements on Cartesian meshes,
therefore there are restrictions on the geometry of the domain and
a higher computational effort compared to $P_k$-type elements.
However, for $P_k$ elements there are indications that such dG methods
do not achieve the optimal order of accuracy,
and it is also known that some of the required properties of the $L^2$-projectors
(e.g.\ property (A2) in \cite[Lemma A1]{Lyu:21}) do not generally hold
in the multidimensional case, especially not for non-tensor product meshes \cite{Oswald:10}.

Let $\Omega:=(r,s)\times(p,q)$, $r<s$, $p<q$, be a rectangular domain in $\R^2$
with boundary $\Gamma$ and unit outward normal $\mathbf{n}$.
As usual,
$\DD = \DD(\xx, t)$, $\BB = \BB(\xx, t)$, $\EE = \EE(\xx, t)$ and $\HH = \HH(\xx, t)$
represent the electric displacement field, magnetic induction, electric and
magnetic field intensities, resp., where $\xx\in\Omega$ and
the time variable $t$ ranges in some interval $(0, T)$, $T > 0$.
Given an electric current density $\JJ = \JJ(\xx, t)$,
we write the transient Maxwell's equations as
\begin{align}
\partial_{t}\DD-\nabla\times \HH=\JJ \quad\text{in }\Omega\times(0,T),\label{non:1_dg}\\
\mu_{0}\partial_{t}\HH+\nabla\times\EE = 0 \quad\text{in }\Omega\times(0,T),\label{non:2_dg}
\end{align}
where
\[
\DD=\eps_{0}\Big( (1+\chi^{(1)})\EE +\chi^{(3)}\vert\EE\vert^{2}\,\EE\Big).
\]
Here $\eps_{0}>0$ denotes the constant vacuum permittivity,
$\mu_{0}:\,\Omega\to (0,\infty)$ is the permeability,
and $\chi^{(1)},\chi^{(3)}:\;\Omega\to (0,\infty)$ are the media susceptibility coefficients.
We assume that the coefficient functions are bounded almost everywhere,
i.e.\ $\mu_{0},\chi^{(1)},\chi^{(3)}\in L_\infty(\Omega)$.

We will consider the derivative
\begin{equation}\label{eq:D_dg}
\partial_{t}\DD
=\eps_{0}\Big( (1+\chi^{(1)})\partial_{t}\EE
+\chi^{(3)}[\vert\EE\vert^{2}+2\EE\EE^{T}]\partial_{t}\EE\Big).
\end{equation}
as an additional equation.

A perfect electric conductor (PEC) boundary condition on $\Gamma$ is assumed, that is
\begin{equation}\label{eq:PEC}
\mathbf{n}\times\EE=0 \quad\text{on }\Gamma\times(0,T).
\end{equation}
In addition, initial conditions have to be specified:
\[
\EE(\xx,0)=\EE_{0}(\xx) \quad\text{and}\quad \HH(\xx,0)=\HH_0(\xx)
\quad \text{for all }\xx\in\Omega,
\]
where $\EE_{0},\HH_0:\;\Omega\to\R^2$ are given functions, and $\HH_0$ satisfies
\begin{align}
\nabla\cdot(\mu_{0}\HH_0)=0 \quad\text{in }\Omega,
\quad\HH_0\cdot\mathbf{n}=0 \quad\text{on }\Gamma.\label{1.7dg}
\end{align}
The divergence-free condition in \eqref{1.7dg} together with \eqref{non:2_dg} implies that
\[
\nabla\cdot(\mu_{0}\HH)=0 \quad\text{in }\Omega\times(0,T).
\]
Here the Transverse Electric mode is considered, where -- for simplicity --
the direction of propagation coincides with the direction of the $z$-axis,
i.e.\ essentially we deal with a two-dimensional problem in space (TE$_z$-mode).
This restriction is only used to simplify the presentation technically;
analogous results for the three-dimensional case can be obtained, too.
The fields reduce to
$\DD=(D_{x},D_{y})$, $\EE=(E_{x},E_{y})$,
$\nabla\times\EE=\partial_{x}E_{y}-\partial_{y}E_{x}$,
$\HH=H_{z}$, $\nabla\times\HH=(\partial_{y}H_{z},-\partial_{x}H_{z})^{T}$,
and $\JJ=(J_{x},J_{y})$,
where the subscripts $x$, $y$ and $z$ denote the $x$-component, $y$-component,
and $z$-component of the vector field, respectively.
In addition we write $\xx:=(x,y)^T$ and $\vert\EE\vert^{2}:=E_{x}^{2}+E_{y}^{2}$.
Then
\[
\EE\EE^{T}\partial_{t}\EE
=\begin{pmatrix}
E_{x}^{2}\partial_{t}E_{x}+E_{x}E_{y}\partial_{t}E_{y} \\
E_{x}E_{y}\partial_{t}E_{x}+E_{y}^{2}\partial_{t}E_{y}
\end{pmatrix},
\]
and the equations \eqref{non:1_dg}--\eqref{eq:D_dg} take the form
\begin{equation}\label{eq:2dsystem}
\begin{aligned}
\partial_{t}D_{x}&=\partial_{y}H_{z}+J_{x}, \\
\partial_{t}D_{y}&=-\partial_{x}H_{z}+J_{y}, \\
\mu_{0}\partial_{t}H_{z}&=-\partial_{x}E_{y}+\partial_{y}E_{x},\\
\partial_{t}D_{x}&=\eps_{0}\Big( (1+\chi^{(1)})\partial_{t}E_{x}
+\chi^{(3)}[\vert\EE\vert^2\partial_{t}E_{x}\\
&\quad +2\big(E_{x}^{2}\partial_{t}E_{x}
+E_{x}E_{y}\partial_{t}E_{y}\big)]\Big), \\
\partial_{t}D_{y}&=\eps_{0}\Big( (1+\chi^{(1)})\partial_{t}E_{y}
+\chi^{(3)}[\vert\EE\vert^2\partial_{t}E_{y}\\
&\quad +2\big(E_{x}E_{y}\partial_{t}E_{x}
+E_{y}^{2}\partial_{t}E_{y}\big)]\Big). \\
\end{aligned}
\end{equation}
The corresponding initial conditions are
\[
E_{x}(\xx,0)=E_{x}^{0}(\xx),
\quad
E_{y}(\xx,0)=E_{y}^{0}(\xx)
\quad\text{and}\quad
H_{z}(\xx,0)=E_{z}^{0}(\xx).
\]
The PEC condition \eqref{eq:PEC} reads as
\begin{equation}\label{eq:pec_con}
E_{x}(\xx,t)\vert_{y=p,q} = E_{y}(\xx,t)\vert_{x=r,s}=0.
\end{equation}

\section{The Nonlinear Electromagnetic Energy at the Continuous Level}
According to the particular structure of the nonlinearity,
a ``nonlinear'' electromagnetic energy of the system \eqref{eq:2dsystem}
can be defined by
\[
\mathcal{E}(t):=\|\EE(t)\|_{\eps_{0}(1+\chi^{(1)})}^2+\|H_{z}(t)\|_{\mu_{0}}^{2}
+\frac{3}{2}\big\|\vert\EE(t)\vert^2\big\|_{\eps_{0}\chi^{(3)}}^{2},
\]
$t\in [0,T)$, where we have used the notation
\[
\|\EE(t)\|_\omega:=\bigg(\int_\Omega \vert\EE(\xx,t)\vert^2\omega(\xx) d\xx \bigg)^{1/2}
\]
for a given weight function $\omega:\;\Omega\to (0,\infty)$.
In the case $\omega= 1,$ the subscript is omitted.

The following theorem demonstrates that the nonlinear electromagnetic energy
is a conservative quantity.
\begin{theorem}\label{th:dg_energycontin}
If $(E_{x},E_{y},H_{z})^T$ is the weak solution
of the system \eqref{eq:2dsystem} in the case of no sources, i.e.\ $\JJ=0$,
then the nonlinear electromagnetic energy of the system \eqref{eq:2dsystem}
at any time $t\in [0,T)$ satisfies
\[
\mathcal{E}(t)=\mathcal{E}(0) = \|\EE_0\|_{\eps_{0}(1+\chi^{(1)})}^2
+\|H_{z}^{0}\|_{\mu_{0}}^{2}+\frac{3}{2}\big\|\vert\EE_0\vert^2\big\|_{\eps_{0}\chi^{(3)}}^{2}.
\]
\end{theorem}
We skip the proof since the details are similar to the more
complicated proof of the semi-discrete energy law (Thm.~\ref{th:dg_weaKenergy}).

\medskip
The domain $\Omega$ is partitioned into rectangular cells $K_{ij}:=I_{i}\times J_{j}$ with
$I_{i}:=(x_{i-\frac{1}{2}},x_{i+\frac{1}{2}})$, $i=1,2,3,\ldots , N_{x}$, and
$J_{j}:=(y_{j-\frac{1}{2}},y_{j+\frac{1}{2}})$, $j=1,2,3,\ldots , N_{y}$,
where
\[
r=:x_{\frac{1}{2}}<x_{\frac{3}{2}}<\ldots <x_{N_{x}+\frac{1}{2}}:=s,
\]\[
p=:y_{\frac{1}{2}}<y_{\frac{3}{2}}<\ldots <y_{N_{y}+\frac{1}{2}}:=q.
\]
The mesh sizes are denoted by $h_{i}^{x}:=x_{i+\frac{1}{2}}-x_{i-\frac{1}{2}}$ and
$h_{j}^{y}:=y_{j+\frac{1}{2}}-y_{i-\frac{1}{2}}$ with
$h_{x}^{\max}:=\max_{1\le i\le N_{x}}h_{i}^{x}$ and $h_{y}^{\max}:=\max_{1\le j\le N_{y}}h_{j}^{y}$.
The maximal mesh size is defined by
\linebreak
$h:=\max\{h_{x}^{\max}, h_{y}^{\max}\}$.
We assume that the mesh is shape-regular,
i.e., if $\varrho_{K_{ij}}$ denotes the radius of the biggest circle contained in $K_{ij}$, we have
$h_{i}^{x}h_{j}^{y}\le C_{sr}\varrho_{K_{ij}}$ for all $K_{ij}$ with a positive constant $C_{sr}$.
The family of cells is denoted by
$\Th:=\{K_{ij}\}_{\substack{i=1,2,3,\ldots, N_{x}\\ j=1,2,3,\ldots , N_{y}}}$.

The finite element space $U_{h}^{k}$ is the space of tensor products of
piecewise polynomials of degree at most $k\in\N$ in each variable on every element $K_{ij}$:
\[
U_{h}^{k}:=\{u:\;u\vert_K\in Q_k(K) \quad\text{for all } K\in\Th\},
\]
where the local space $Q_k(K)$ consists of tensor products of univariate polynomials
of degree up to $k$ on a cell $K$.
Note that $U_{h}^{k}\not\subset C(\overline{\Omega})$ in general.

The numerical approximation of a function $u:\;\overline{\Omega}\to\R$ is denoted
by $u_{h}\in U_h^k$.
The limiting value of $u_h$ at $x_{i+\frac{1}{2}}$ from the right cell $K_{i+1,j}$
is denoted by $u_{h}(x_{i+\frac{1}{2}}^{+}, y)$,
$(u_{h})_{i+\frac{1}{2},y}^{+}$ or $u_{h}^{+}(x_{i+\frac{1}{2}}, y)$, and from
the left cell $K_{ij}$ by $u_{h}(x_{i+\frac{1}{2}}^{-}, y)$,
$(u_{h})_{i+\frac{1}{2},y}^{-}$ or $u_{h}^{-}(x_{i+\frac{1}{2}}, y)$,
respectively.
An analogous convention is used in the $y$-direction.

The numerical fluxes are obtained by means of integration by parts.
They should be considered and designed carefully to ensure conservation of energy,
numerical stability and optimal accuracy of the approximate solution.
The numerical flux densities are functions that are defined on the cell boundaries.
The alternating flux densities are defined in a simple and elegant way like in
LDG (local discontinuous Galerkin) methods for diffusion equations,
second order wave equation and Maxwell's equations \cite{Cockburn:98,Xing:13,Li:17b}.
Fixing a constant $c_{0}>0$ independent of $h$, the alternating flux densities are:
\begin{align}
\hat{E}_{x,h}(x,y_{j+\frac{1}{2}})&:=E_{x,h}^{+}(x,y_{j+\frac{1}{2}})
\quad\text{for all } j = 1, 2, 3,\ldots,N_{y}-1,\nonumber\\
\hat{E}_{x,h}(x,y_{\frac{1}{2}})&:=\hat{E}_{x,h}(x,y_{N_{y}+\frac{1}{2}}):=0,\nonumber\\
\hat{E}_{y,h}(x_{i+\frac{1}{2}},y)&:=E_{y,h}^{+}(x_{i+\frac{1}{2}},y)
\quad\text{for all } i= 1, 2, 3,\ldots , N_{x}-1,\nonumber\\
\hat{E}_{y,h}(x_{\frac{1}{2}},y)&:=\hat{E}_{y,h}(x_{N_{x}+\frac{1}{2}},y):=0,\nonumber\\
\hat{H}_{z,h}(x,y_{j+\frac{1}{2}})&:=H_{z,h}^{-}(x,y_{j+\frac{1}{2}})
\quad\text{for all } j = 1, 2, 3,\ldots , N_{y},\nonumber\\
\hat{H}_{z,h}(x,y_{\frac{1}{2}})&:=H_{z,h}^{+}(x,y_{\frac{1}{2}})
+c_{0}\llbracket E_{x,h}(x,y_{\frac{1}{2}}) \rrbracket,\label{eq:jump6}\\
\hat{H}_{z,h}(x_{i+\frac{1}{2}},y)&:=H_{z,h}^{-}(x_{i+\frac{1}{2}},y)
\quad\text{for all } i = 1, 2, 3,\ldots, N_{x},\nonumber\\
\hat{H}_{z,h}(x_{\frac{1}{2}},y)&:=H_{z,h}^{+}(x_{\frac{1}{2}},y)
-c_{0}\llbracket E_{y,h}(x_{\frac{1}{2}},y) \rrbracket,\label{eq:jump8}
\end{align}
where the jump terms in \eqref{eq:jump6}, \eqref{eq:jump8} are defined as
\[
\llbracket E_{x,h}(x,y_{\frac{1}{2}}) \rrbracket :=E_{x,h}^{+}(x,y_{\frac{1}{2}})-0,
\qquad
\llbracket E_{y,h}(x_{\frac{1}{2}},y) \rrbracket :=E_{y,h}^{+}(x_{\frac{1}{2}},y)-0.
\]
On interior cell boundaries, the jumps are denoted by
\[
\llbracket \Psi \rrbracket:=\Psi^{+}-\Psi^{-}.
\]
For $c_{0}=\frac{1}{2}$, the flux densities \eqref{eq:jump6}, \eqref{eq:jump8} match with
the standard upwind flux densities
\begin{align*}
&\hat{H}_{z,h}(x,y_{\frac{1}{2}})
:=\frac{1}{2}[H_{z,h}^{+}(x,y_{\frac{1}{2}})
+H_{z,h}^{-}(x,y_{\frac{1}{2}})]
+\frac{1}{2}\llbracket E_{x,h}(x,y_{\frac{1}{2}}) \rrbracket,\\
&\hat{H}_{z,h}(x_{\frac{1}{2}},y)
:=\frac{1}{2}[H_{z,h}^{+}(x_{\frac{1}{2}},y)
+H_{z,h}^{-}(x_{\frac{1}{2}},y)]
-\frac{1}{2}\llbracket E_{y,h}(x_{\frac{1}{2}},y) \rrbracket,
\end{align*}
where the undefined $H_{z,h}^{-}(x,y_{\frac{1}{2}})$ and $H_{z,h}^{-}(x_{\frac{1}{2}},y)$
are replaced by $H_{z,h}^{+}(x,y_{\frac{1}{2}})$ and $H_{z,h}^{+}(x_{\frac{1}{2}},y)$, respectively.

\bigskip
\section{Spatial Discretization for Discontinuous Galerkin Method}
\label{sec:spatial_discr}
For the test functions $(\Phi_{1h},\Phi_{2h},\Phi_{3h})^T\in (U_{h}^{k})^3$,
the discontinuous Galerkin formulation for the equations \eqref{eq:2dsystem}
with respect to the semi-discrete solution
$(E_{x,h},E_{y,h},H_{z,h})^T\in C^{1}(0,T,U_{h}^{k})^3$
reads as follows
(for shortness, we omit the formal differentials $d\xx$ in the double integrals):
\begin{align}
&\int_{K_{ij}}\partial_t D_{x,h}\,\Phi_{1h}
-\int_{I_{i}}[(\hat{H}_{z,h}\Phi_{1h}^{-})_{x,j+\frac{1}{2}}-(\hat{H}_{z,h}\Phi_{1h}^{+})_{x,j
-\frac{1}{2}}]dx\nonumber\\
&+\int_{K_{ij}}H_{z,h}\partial_{y}\Phi_{1h}
-\int_{K_{ij}}J_{x,h}\Phi_{1h}=0, \label{eq:w_x}\\
&\int_{K_{ij}}\partial_t D_{y,h}\,\Phi_{2h}
+\int_{J_{j}}[(\hat{H}_{z,h}\Phi_{2h}^{-})_{i+\frac{1}{2},y}-(\hat{H}_{z,h}\Phi_{2h}^{+})_{i
-\frac{1}{2},y}]dy\nonumber\\
&-\int_{K_{ij}}H_{z,h}\partial_{x}\Phi_{2h}
-\int_{K_{ij}}J_{y,h}\Phi_{2h}=0, \label{eq:w_y}\\
&\int_{K_{ij}}\mu_{0}\partial_t H_{z,h}\,\Phi_{3h}\nonumber\\
&+\int_{J_{j}}[(\hat{E}_{y,h}\Phi_{3h}^{-})_{i+\frac{1}{2},y}-(\hat{E}_{y,h}\Phi_{3h}^{+})_{i
-\frac{1}{2},y}]dy\nonumber\\
&-\int_{K_{ij}}E_{y,h}\partial_{x}\Phi_{3h}
-\int_{I_{i}}[(\hat{E}_{x,h}\Phi_{3h}^{-})_{x,j+\frac{1}{2}}
-(\hat{E}_{x,h}\Phi_{3h}^{+})_{x,j-\frac{1}{2}}]dx\nonumber\\
&+\int_{K_{ij}}E_{x,h}\partial_{y}\Phi_{3h}=0,\label{eq:w_h}
\end{align}
\begin{align}
\int_{K_{ij}}\partial_t D_{x,h}\,\Phi_{1h}
&=\int_{K_{ij}}\eps_{0}(1+\chi^{(1)})\partial_t E_{x,h}\,\Phi_{1h}
+\int_{K_{ij}}\eps_{0}\chi^{(3)}\Big[\vert\EE_h\vert^2\partial_t E_{x,h}\,\Phi_{1h}
\nonumber\\
&
+2\big(E_{x,h}^{2}\partial_t E_{x,h}\,\Phi_{1h}
+E_{x,h}E_{y,h}\partial_t E_{y,h}\,\Phi_{1h}\big)\Big],\label{eq:w_d1}\\
\int_{K_{ij}}\partial_t D_{y,h}\,\Phi_{2h}
&=\int_{K_{ij}}\eps_{0}(1+\chi^{(1)})\partial_t E_{y,h}\,\Phi_{2h}
+\int_{K_{ij}}\eps_{0}\chi^{(3)}\Big[\vert\EE_h\vert^2\partial_t E_{y,h}\,\Phi_{2h}
\nonumber\\
&
+2\big(E_{y,h}^{2}\partial_t E_{y,h}\,\Phi_{2h}
+E_{x,h}E_{y,h}\partial_t E_{x,h}\,\Phi_{2h}\big)\Big].\label{eq:w_d2}
\end{align}
The initial conditions are defined as
\[
E_{x,h}(\xx,0)=E_{x,h}^{0}(\xx),
\quad
E_{y,h}(\xx,0)=E_{y,h}^{0}(\xx)
\quad\text{and}\quad
H_{z,h}(\xx,0)=E_{z,h}^{0}(\xx),
\]
where the concrete choice of the discrete initial data
$(E_{x,h}^{0},E_{y,h}^{0},H_{z,h}^{0})^T\in (U_{h}^{k})^3$ will be given later.

\bigskip
\section{The Nonlinear Electromagnetic Energy of
the Semi-Discrete Discontinuous Galerkin Discretization}
\label{sec:nonlinear_energy_discr}
The nonlinear electromagnetic energy of the semi-discrete discontinuous Galerkin discretization
of the system \eqref{eq:w_x}--\eqref{eq:w_d2} is defined by
\begin{align*}
\mathcal{E}_{h}(t)
&:=\|\EE_h(t)\|_{\eps_{0}(1+\chi^{(1)})}^2+\|H_{z,h}(t)\|_{\mu_{0}}^{2}
+\frac{3}{2}\big\|\vert\EE_h(t)\vert^2\big\|_{\eps_{0}\chi^{(3)}}^{2}\\
&\qquad +\frac{c_{0}}{2}\int_{0}^{t}\Big[\int_{r}^{s}(E_{x,h}^{+}(t))_{x,\frac{1}{2}}^{2}dx
+\int_{p}^{q}(E_{y,h}^{+}(t))_{\frac{1}{2},y}^{2}dy\Big],
\end{align*}
$t\in [0,T)$.
In the next, we will show that the nonlinear electromagnetic energy at the semi-discrete level
of the system \eqref{eq:w_x}--\eqref{eq:w_d2} at time $t$ is conserved and bounded.
\begin{theorem}\label{th:dg_weaKenergy}
Let $(E_{x,h},E_{y,h},H_{z,h})^T\in C^{1}(0,T,U_{h}^{k})^3$
be the semi-discrete solution of the system \eqref{eq:w_x}--\eqref{eq:w_d2}
for given $\JJ_h\in C(0,T,U_{h}^{k})^2$,
then the nonlinear electromagnetic energy of the system \eqref{eq:w_x}--\eqref{eq:w_d2}
for the vanishing current density at any time $t\in [0,T)$ satisfies
\begin{equation}\label{eq:semi_discrete_zero}
\mathcal{E}_{h}(t)=\mathcal{E}_{h}(0)
=\|\EE_h^0\|_{\eps_{0}(1+\chi^{(1)})}^{2}+\|H_{z,h}^{0}\|_{\mu_0}^{2}
+\frac{3}{2}\big\|\vert\EE_h^0\vert^2\big\|_{\eps_{0}\chi^{(3)}}^{2},
\end{equation}
and for non-zero current density
\begin{equation}\label{eq:semi_discrete_nonzero}
\mathcal{E}_h(t)\le 2\mathcal{E}_h(0)
+ 8\bigg(\int_{0}^{t}\|\JJ_h(s)\|_{(\eps_{0}(1+\chi^{(1)}))^{-1}}ds\bigg)^2.
\end{equation}
\end{theorem}
\proof
Taking $ \Phi_{1h}:= E_{x,h}$ in the equations \eqref{eq:w_x} and \eqref{eq:w_d1}, and substituting
the equation \eqref{eq:w_d1} into the equation \eqref{eq:w_x}, we have
\begin{equation}\label{eq:w_1}
\begin{aligned}
&\int_{K_{ij}}\eps_{0}(1+\chi^{(1)})\partial_t E_{x,h}\,E_{x,h}
+\int_{K_{ij}}\eps_{0}\chi^{(3)}\Big[\vert\EE_h\vert^2\partial_t E_{x,h}\,E_{x,h}\\
&\quad +2\big(E_{x,h}^{2}\partial_t E_{x,h}\,E_{x,h}+E_{x,h}E_{y,h}\partial_t E_{y,h}\,E_{x,h}\big)\Big]\\
&-\int_{I_{i}}[(\hat{H}_{z,h}E_{x,h}^{-})_{x,j+\frac{1}{2}}
-(\hat{H}_{z,h}E_{x,h}^{+})_{x,j-\frac{1}{2}}]dx
+\int_{K_{ij}}H_{z,h}\partial_{y}E_{x,h}-\int_{K_{ij}}J_{x,h}E_{x,h}=0.
\end{aligned}
\end{equation}
Taking $ \Phi_{2h}:= E_{y,h}$ in the equations \eqref{eq:w_y} and \eqref{eq:w_d2}, and substituting
the equation \eqref{eq:w_d2} into the equation \eqref{eq:w_y}, we have
\begin{equation}\label{eq:w_2}
\begin{aligned}
&\int_{K_{ij}}\eps_{0}(1+\chi^{(1)})\partial_t E_{y,h}\,E_{y,h}
+\int_{K_{ij}}\eps_{0}\chi^{(3)}\Big[\vert\EE_h\vert^2\partial_t E_{y,h}\,E_{y,h}\\
&\quad +2\big(E_{y,h}^{2}\partial_t E_{y,h}\,E_{y,h}+E_{x,h}E_{y,h}\partial_t E_{x,h}\,E_{y,h}\big)\Big]\\
&+\int_{J_{j}}[(\hat{H}_{z,h}E_{y,h}^{-})_{i+\frac{1}{2},y}
-(\hat{H}_{z,h}E_{y,h}^{+})_{i-\frac{1}{2},y}]dy
-\int_{K_{ij}}H_{z,h}\partial_{x}E_{y,h}-\int_{K_{ij}}J_{y,h}E_{y,h} =0.
\end{aligned}
\end{equation}
Adding the equations \eqref{eq:w_1} and \eqref{eq:w_2}, we get
\begin{align*}
&\frac{1}{2}\frac{d}{dt}\int_{K_{ij}}\eps_{0}(1+\chi^{(1)})\vert\EE_h\vert^2
+\int_{K_{ij}}\eps_{0}\chi^{(3)}\vert\EE_h\vert^2\frac{1}{2}\partial_t \vert\EE_h\vert^2\\
&+\int_{K_{ij}}\eps_{0}\chi^{(3)}2\Big[E_{x,h}^{2}\frac{1}{2}\partial_t E_{x,h}^{2}
+E_{y,h}^{2}\frac{1}{2}\partial_t E_{y,h}^{2}\Big]\\
&+\int_{K_{ij}}\eps_{0}\chi^{(3)}2\Big[E_{x,h}E_{y,h}[\partial_t E_{y,h}\,E_{x,h}
+\partial_t E_{x,h}\,E_{y,h}]\Big]\\
&-\int_{I_{i}}[(\hat{H}_{z,h}E_{x,h}^{-})_{x,j+\frac{1}{2}}
-(\hat{H}_{z,h}E_{x,h}^{+})_{x,j-\frac{1}{2}}]dx\\
&+\int_{J_{j}}[(\hat{H}_{z,h}E_{y,h}^{-})_{i+\frac{1}{2},y}
-(\hat{H}_{z,h}E_{y,h}^{+})_{i-\frac{1}{2},y}]dy\\
&+\int_{K_{ij}}H_{z,h}\partial_{y}E_{x,h}-\int_{K_{ij}}H_{z,h}\partial_{x}E_{y,h}\\
&-\int_{K_{ij}}J_{x,h}E_{x,h}-\int_{K_{ij}}J_{y,h}E_{y,h}=0.
\end{align*}
The integrands corresponding to the cubic nonlinearity can be rewritten as follows:
\begin{align*}
&\vert\EE_h\vert^2\frac{1}{2}\partial_t \vert\EE_h\vert^2
+2\Big[E_{x,h}^{2}\frac{1}{2}\partial_t E_{x,h}^{2}+E_{y,h}^{2}\frac{1}{2}\partial_t E_{y,h}^{2}\Big]\\
&\quad +2\Big[E_{x,h}E_{y,h}[\partial_t E_{y,h}\,E_{x,h}+\partial_t E_{x,h}\,E_{y,h}]\Big]\\
&=\frac{1}{4}\partial_t \vert\EE_h\vert^4
+\frac{1}{2}\Big[\partial_t E_{x,h}^4+\partial_t E_{y,h}^4\Big]\\
&\quad +2\Big[E_{x,h}E_{y,h}\partial_t (E_{y,h}E_{x,h})]\Big]\\
&=\frac{1}{4}\partial_t \vert\EE_h\vert^4
+\frac{1}{2}\partial_t (E_{x,h}^4+ E_{y,h}^4)
+\partial_t (E_{y,h}E_{x,h})^2\\
&=\frac{1}{4}\partial_t \vert\EE_h\vert^4+\frac{1}{2}\partial_t \vert\EE_h\vert^4
=\frac{3}{4}\partial_t \vert\EE_h\vert^4.
\end{align*}
Thus we arrive at
\begin{equation}\label{eq:w_3}
\begin{aligned}
&\frac{1}{2}\frac{d}{dt}\int_{K_{ij}}\eps_{0}(1+\chi^{(1)})\vert\EE_h\vert^2
+\frac{3}{4}\frac{d}{dt}\int_{K_{ij}} \eps_{0}\chi^{(3)}\vert\EE_h\vert^4\\
&-\int_{I_{i}}[(\hat{H}_{z,h}E_{x,h}^{-})_{x,j+\frac{1}{2}}
-(\hat{H}_{z,h}E_{x,h}^{+})_{x,j-\frac{1}{2}}]dx\\
&+\int_{J_{j}}[(\hat{H}_{z,h}E_{y,h}^{-})_{i+\frac{1}{2},y}
-(\hat{H}_{z,h}E_{y,h}^{+})_{i-\frac{1}{2},y}]dy\\
&+\int_{K_{ij}}H_{z,h}\partial_{y}E_{x,h}-\int_{K_{ij}}H_{z,h}\partial_{x}E_{y,h}\\
&-\int_{K_{ij}}J_{x,h}E_{x,h}-\int_{K_{ij}}J_{y,h}E_{y,h}=0.
\end{aligned}
\end{equation}
Taking $\Phi_{3h}:= H_{z,h}$ in the equations \eqref{eq:w_h}, we have
\begin{equation}\label{eq:w_4}
\begin{aligned}
&\int_{K_{ij}}\mu_{0}\partial_t H_{z,h}\,H_{z,h}
+\int_{J_{j}}[(\hat{E}_{y,h}H_{z,h}^{-})_{i+\frac{1}{2},y}\\
&-(\hat{E}_{y,h}H_{z,h}^{+})_{i-\frac{1}{2},y}]dy
-\int_{K_{ij}}E_{y,h}\partial_{x}H_{z,h}\\
&-\int_{I_{i}}[(\hat{E}_{x,h}H_{z,h}^{-})_{x,j+\frac{1}{2}}
-(\hat{E}_{x,h}H_{z,h}^{+})_{x,j-\frac{1}{2}}]dx
+\int_{K_{ij}}E_{x,h}\partial_{y}H_{z,h}=0.
\end{aligned}
\end{equation}
Adding the equations \eqref{eq:w_3} and \eqref{eq:w_4}, we obtain
\begin{equation}\label{eq:w_7}
\begin{aligned}
&\frac{1}{2}\frac{d}{dt}\int_{K_{ij}}\eps_{0}(1+\chi^{(1)})\vert\EE_h\vert^2
+\frac{1}{2}\frac{d}{dt}\int_{K_{ij}}\mu_{0}H_{z,h}^{2}
+\frac{3}{4}\frac{d}{dt}\int_{K_{ij}} \eps_{0}\chi^{(3)}\vert\EE_h\vert^4\\
&-\int_{I_{i}}[(\hat{H}_{z,h}E_{x,h}^{-})_{x,j+\frac{1}{2}}
-(\hat{H}_{z,h}E_{x,h}^{+})_{x,j-\frac{1}{2}}]dx\\
&+\int_{J_{j}}[(\hat{H}_{z,h}E_{y,h}^{-})_{i+\frac{1}{2},y}
-(\hat{H}_{z,h}E_{y,h}^{+})_{i-\frac{1}{2},y}]dy\\
&+\int_{J_{j}}[(\hat{E}_{y,h}H_{z,h}^{-})_{i+\frac{1}{2},y}
-(\hat{E}_{y,h}H_{z,h}^{+})_{i-\frac{1}{2},y}]dy\\
&-\int_{I_{i}}[(\hat{E}_{x,h}H_{z,h}^{-})_{x,j+\frac{1}{2}}
-(\hat{E}_{x,h}H_{z,h}^{+})_{x,j-\frac{1}{2}}]dx\\
&+\int_{K_{ij}}H_{z,h}\partial_{y}E_{x,h}-\int_{K_{ij}}H_{z,h}\partial_{x}E_{y,h}\\
&-\int_{K_{ij}}E_{y,h}\partial_{x}H_{z,h}+\int_{K_{ij}}E_{x,h}\partial_{y}H_{z,h}\\
&-\int_{K_{ij}}J_{x,h}E_{x,h}-\int_{K_{ij}}J_{y,h}E_{y,h}=0.
\end{aligned}
\end{equation}
In the next step, the equations \eqref{eq:w_7} are summed up with respect to
both indices $1\le i\le N_{x}$ and $1\le j\le N_{y}$.
The sums resulting from the terms on the second to fourth lines allow
the following simplification, see \cite[eqs.~(3.18)--(3.19)]{Li:17b}:
\begin{equation}\label{eq:identities12}
\begin{aligned}
&\sum_{j=1}^{N_{y}}\Big[-\int_{I_{i}}[(\hat{H}_{z,h}E_{x,h}^{-})_{x,j+\frac{1}{2}}
-(\hat{H}_{z,h}E_{x,h}^{+})_{x,j-\frac{1}{2}}]dx\\
&-\int_{I_{i}}[(\hat{E}_{x,h}H_{z,h}^{-})_{x,j+\frac{1}{2}}
-(\hat{E}_{x,h}H_{z,h}^{+})_{x,j-\frac{1}{2}}]dx\\
&+\int_{K_{ij}}H_{z,h}\partial_{y}E_{x,h}+\int_{K_{ij}}E_{x,h}\partial_{y}H_{z,h}\Big]
=c_{0}\int_{I_{i}}(E_{x,h}^{+})_{x,\frac{1}{2}}^{2}dx,\\
&\sum_{i=1}^{N_{x}}\Big[\int_{J_{j}}[(\hat{H}_{z,h}E_{y,h}^{-})_{i+\frac{1}{2},y}
-(\hat{H}_{z,h}E_{y,h}^{+})_{i-\frac{1}{2},y}]dy\\
&+\int_{J_{j}}[(\hat{E}_{y,h}H_{z,h}^{-})_{i+\frac{1}{2},y}
-(\hat{E}_{y,h}H_{z,h}^{+})_{i-\frac{1}{2},y}]dy\\
&-\int_{K_{ij}}H_{z,h}\partial_{x}E_{y,h}-\int_{K_{ij}}E_{y,h}\partial_{x}H_{z,h}\Big]
=c_{0}\int_{J_{j}}(E_{y,h}^{+})_{\frac{1}{2},y}^{2}dy.
\end{aligned}
\end{equation}
Using these relationships, we get
\begin{equation}\label{eq:w_8}
\begin{aligned}
&\frac{1}{2}\frac{d}{dt}\|\EE_h\|_{\eps_{0}(1+\chi^{(1)})}^2
+\frac{1}{2}\frac{d}{dt}\|H_{z,h}\|_{\mu_{0}}^{2}
+\frac{3}{4}\frac{d}{dt}\big\|\vert\EE_h\vert^2\big\|_{\eps_{0}\chi^{(3)}}^{2}\\
&+c_{0}\int_{r}^{s}(E_{x,h}^{+})_{x,\frac{1}{2}}^{2}dx
+c_{0}\int_{p}^{q}(E_{y,h}^{+})_{\frac{1}{2},y}^{2}dy
=\int_{\Omega}[J_{x,h}E_{x,h}+J_{y,h}E_{y,h}].
\end{aligned}
\end{equation}
The right-hand side of the equation \eqref{eq:w_8} is estimated by means of
Cauchy-Schwarz inequalities as follows:
\begin{align*}
&\int_{\Omega}[J_{x,h}E_{x,h}+J_{y,h}E_{y,h}]
\le \int_{\Omega} \vert\JJ_h\vert\vert\EE_h\vert\\
&= \int_{\Omega} \vert\EE_h\vert\sqrt{\eps_{0}(1+\chi^{(1)})}\vert\JJ_h\vert
\sqrt{(\eps_{0}(1+\chi^{(1)}))^{-1}}\\
&\le \|\EE_h\|_{\eps_{0}(1+\chi^{(1)})}\|\JJ_h\|_{(\eps_{0}(1+\chi^{(1)}))^{-1}}.
\end{align*}
Then we obtain from equation \eqref{eq:w_8}
\begin{align*}
&\frac{1}{2}\frac{d}{dt}\|\EE_h\|_{\eps_{0}(1+\chi^{(1)})}^2
+\frac{1}{2}\frac{d}{dt}\|H_{z,h}\|_{\mu_{0}}^{2}
+\frac{3}{4}\frac{d}{dt}\big\|\vert\EE_h\vert^2\big\|_{\eps_{0}\chi^{(3)}}^{2}\\
&+c_{0}\int_{r}^{s}(E_{x,h}^{+})_{x,\frac{1}{2}}^{2}dx
+c_{0}\int_{p}^{q}(E_{y,h}^{+})_{\frac{1}{2},y}^{2}dy\\
&\le \|\EE_h\|_{\eps_{0}(1+\chi^{(1)})}\|\JJ_h\|_{(\eps_{0}(1+\chi^{(1)}))^{-1}}.
\end{align*}
Integration of both sides from $0$ to $t$ yields
\begin{align*}
\frac{1}{2}\mathcal{E}_h(t)-\frac{1}{2}\mathcal{E}_h(0)
&\le \int_{0}^{t}\|\EE_h(s)\|_{\eps_{0}(1+\chi^{(1)})}
\|\JJ_h(s)\|_{(\eps_{0}(1+\chi^{(1)}))^{-1}}ds,
\end{align*}
hence
\[
\mathcal{E}_h(t)\le\mathcal{E}_h(0)
+ 2\int_{0}^{t}\sqrt{\mathcal{E}_h(s)}\|\JJ_h(s)\|_{(\eps_{0}(1+\chi^{(1)}))^{-1}}ds.
\]
Then the Gronwall--Ou-Iang's inequality \cite{Pachpatte:94} implies that
\[
\sqrt{\mathcal{E}_h(t)}\le\sqrt{\mathcal{E}_h(0)}
+ 2\int_{0}^{t}\|\JJ_h(s)\|_{(\eps_{0}(1+\chi^{(1)}))^{-1}}ds.
\]
Squaring this estimate together with an elementary inequality
completes the proof of \eqref{eq:semi_discrete_nonzero}.
The relationship \eqref{eq:semi_discrete_zero} immediately follows
from integration of \eqref{eq:w_8} for the case $\JJ_h=0$.
\close

\bigskip
\section{Error Estimates for the Semi-Discrete Discontinuous Galerkin Discretization}
Projection operators play an important role in the error analysis,
and we will begin with defining 1D projectors that are frequently used
in discontinuous Galerkin methods \cite{Cockburn:01c,Dong:09}.
In this presentation, we closely follow \cite{Li:17b}.
Let $\mathcal{P}_k(I_{i})$ denote the $k$th degree polynomial space over the interval $I_{i}$,
$k\in\N$.
For any function $u\in H^1(I_{i})$, we define
\[
P_{x}^{\pm}:\; H^1(I_{i})\to \mathcal{P}_k(I_{i})
\]
by
\begin{align*}
&\int_{I_{i}} (P_{x}^{+} u)w dx = \int_{I_{i}} u\,w dx\\
&\quad\text{for all } w\in P_{k-1}(I_{i})\quad\text{and}\quad
P_{x}^{+}u\Big(x_{i-\frac{1}{2}}^{+}\Big) := u\Big(x_{i-\frac{1}{2}}^{+}\Big),
\\
&\int_{I_{i}} (P_{x}^{-} u)w dx = \int_{I_{i}} u\,w dx\\
&\quad\text{for all } w\in P_{k-1}(I_{i})
\quad\text{and}\quad
P_{x}^{-}u\Big(x_{i+\frac{1}{2}}^{-}\Big) := u\Big(x_{i+\frac{1}{2}}^{-}\Big).
\end{align*}
Analogously, for any function $u\in H^1(J_{j})$, the projection operators in $y$-direction
\[
P_{y}^{\pm}:\; H^1(J_{j})\to P_{k}(J_{j})
\]
are defined by
\begin{align*}
&\int_{J_{j}} (P_{y}^{+} u) w dy = \int_{J_{j}} u\,w dy\\
&\quad\text{for all } w\in P_{k-1}(J_{j})
\quad\text{and}\quad
P_{y}^{+}u\Big(y_{j-\frac{1}{2}}^{+}\Big) := u\Big(y_{j-\frac{1}{2}}^{+}\Big),
\\
&\int_{J_{j}} (P_{y}^{-} u) wdy = \int_{J_{j}} u\,w dy\\
&\quad\text{for all } w\in P_{k-1}(J_{j})
\quad\text{and}\quad
P_{y}^{-}u\Big(y_{j+\frac{1}{2}}^{-}\Big) := u\Big(y_{j+\frac{1}{2}}^{-}\Big).
\end{align*}
The standard local $L_2$-projection operators in 1D are denoted by
\[
P_{x}:\; H^1(I_{i})\to \mathcal{P}_k(I_{i})
\quad\text{and}\quad
P_{y}:\; H^1(J_{j})\to P_{k}(J_{j}).
\]
The 2D projection operators for the rectangular elements $K_{ij}=I_{i}\times J_{j}$
are defined as tensor products of the 1D projectors. We define
\begin{equation}
\Pi_{1}:= P_{x}\otimes P_{y}^{+}:\; H^{2}(K_{ij})\to Q_{k}(K_{ij}),\label{eq:pi_dg1}
\end{equation}
which satisfies
\begin{align*}
\int_{K_{ij}}[ \Pi_{1}w(x,y)\partial_y u_{h}(x,y)]
&=\int_{K_{ij}}[ w(x,y)\partial_y u_{h}(x,y)],\\
\int_{I_{i}}\Pi_{1}w\Big(x,y_{j-\frac{1}{2}}^{+}\Big)u_{h}\Big(x,y_{j-\frac{1}{2}}^{+}\Big)dx
&=\int_{I_{i}}w\Big(x,y_{j-\frac{1}{2}}^{+}\Big)u_{h}\Big(x,y_{j-\frac{1}{2}}^{+}\Big)dx
\end{align*}
for all $w\in H^{2}(K_{ij})$ and $u_{h}\in Q_{k}(K_{ij})$ \cite{Meng:16,Li:17b}.
The projection $\Pi_{2}$ is defined as
\begin{equation}\label{eq:pi_dg2}
\Pi_{2}:= P_{x}^{+}\otimes P_{y}:\; H^{2}(K_{ij})\to Q_{k}(K_{ij})
\end{equation}
and satisfies
\begin{align*}
\int_{K_{ij}}[ \Pi_{2}w(x,y)\partial_x u_{h}(x,y)]
&=\int_{K_{ij}}[ w(x,y)\partial_x u_{h}(x,y)],\\
\int_{J_{j}}\Pi_{2}w\Big(x_{i-\frac{1}{2}}^{+},y\Big)u_{h}\Big(x_{i-\frac{1}{2}}^{+},y\Big)dy
&=\int_{J_{j}}w\Big(x_{i-\frac{1}{2}}^{+},y\Big)u_{h}\Big(x_{i-\frac{1}{2}}^{+},y\Big)dy
\end{align*}
for all $w\in H^{2}(K_{ij})$ and $u_{h}\in Q_{k}(K_{ij})$.
The projection $\Pi_{3}$ is defined as
\begin{equation}\label{eq:pi_dg3}
\Pi_{3}:= P_{x}^{-}\otimes P_{y}^{-}:\; H^{2}(K_{ij})\to Q_{k}(K_{ij}).
\end{equation}
It satisfies
\begin{align*}
\int_{K_{ij}}[ \Pi_{3}w(x,y) u_{h}(x,y)]
&=\int_{K_{ij}}[ w(x,y)u_{h}(x,y)],\\
\int_{I_{i}}\Pi_{3}w\Big(x,y_{j+\frac{1}{2}}^{-}\Big)u_{h}\Big(x,y_{j+\frac{1}{2}}^{-}\Big)dx
&=\int_{I_{i}}w\Big(x,y_{j+\frac{1}{2}}^{-}\Big)u_{h}\Big(x,y_{j+\frac{1}{2}}^{-}\Big)dx,
\\
\int_{J_{j}}\Pi_{3}w\Big(x_{i+\frac{1}{2}}^{-},y\Big)u_{h}\Big(x_{i+\frac{1}{2}}^{-},y\Big)dy
&=\int_{J_{j}}w\Big(x_{i+\frac{1}{2}}^{-},y\Big)u_{h}\Big(x_{i+\frac{1}{2}}^{-},y\Big)dy,\\
\Pi_{3}w\Big(x_{i+\frac{1}{2}}^{-},y_{j+\frac{1}{2}}^{-}\Big)
&=w\Big(x_{i+\frac{1}{2}}^{-},y_{j+\frac{1}{2}}^{-}\Big)
\end{align*}
for all $w\in H^{2}(K_{ij})$ and $u_{h}\in Q_{k-1}(K_{ij})$.
The use of the $H^{2}$-spaces for the point values makes sense due to
the Sobolev embedding $H^{2} \subset C^{0}$ in 2D.
The 2D $L_2$-projection operator is usually defined by
\begin{equation}
\Pi_{4}:= P_{x}\otimes P_{y}:\; H^{2}(K_{ij})\to Q_{k}(K_{ij}),\label{eq:pi_dg4}
\end{equation}
for their properties see \cite{Cockburn:01c,Dong:09}, \cite[eqs.~(3.33)--(3.42)]{Li:17b}.
\begin{lemma}\label{lemma:1d_dg}
If $w$ is a product of 1D functions, i.e.\ $w(x,y)=f(x) g(y)$, where $f\in H^1(I_{i})$
and $g\in H^1(J_{j})$, then
\[
\begin{alignedat}{2}
\Pi_{1}w(x,y)&=P_{x}f(x) P_{y}^{+}g(y), & \Pi_{2}w(x,y)&=P_{x}^{+}f(x) P_{y}g(y),\\
\Pi_{3}w(x,y)&=P_{x}^{-}f(x) P_{y}^{-}g(y),\ & \Pi_{4}w(x,y)&=P_{x}f(x) P_{y}g(y).
\end{alignedat}
\]
\end{lemma}
These results demonstrate that the 2D projections are tensor products of 1D projections,
for details see \cite{Cockburn:01c,Dong:09}.
\begin{lemma}\label{lemma:pro_error_dg}
The projection operators $\Pi_{1},\dots,\Pi_{4}$, defined in
\eqref{eq:pi_dg1}--\eqref{eq:pi_dg4},
have the following property:
For $k\in\N$, there exists a constant $C>0$ independent of $h$ such that
\[ 
\|\Pi_{i}u-u\|\le Ch^{k+1}\|u\|_{H^{k+1}(\Omega)}
\] 
for all $u\in H^{k+1}(\Omega),\ i=1,\ldots,4$.
\end{lemma}
Now we are prepared to derive an error estimate.
Let $(E_{x},E_{y},H_{z})^T$ be the weak solution of \eqref{eq:2dsystem}
and $(E_{x,h},E_{y,h},H_{z,h})^T$ be corresponding numerical solution
of the semi-discrete scheme \eqref{eq:w_x}--\eqref{eq:w_d2}.
We denote the error terms for later use by
\begin{equation}\label{eq:zeta_full_x}
\zeta_{x}:=E_{x}-E_{x,h}
=\eta_{x}-\eta_{x,h},
\end{equation}
where
\begin{equation}\label{eq:etas_full_x}
\eta_{x}:=E_{x}-\Pi_{1} E_{x},
\quad
\eta_{x,h}:= E_{x,h}-\Pi_{1} E_{x}.
\end{equation}
Similarly for the $y$-component of the electric field we set
\begin{equation}\label{eq:zeta_full_y}
\zeta_{y}:=E_{y}-E_{y,h}
=\eta_{y}-\eta_{y,h},
\end{equation}
where
\begin{equation}\label{eq:etas_full_y}
\eta_{y}:=E_{y}-\Pi_{2} E_{y},
\quad
\eta_{y,h}:= E_{y,h}-\Pi_{2} E_{y}.
\end{equation}
The error terms for the magnetic field are defined by:
\begin{equation}\label{eq:xi_full_z}
\xi_{z}:=H_{z}-H_{z,h}
=\theta_{z}-\theta_{z,h},
\end{equation}
where
\begin{equation}\label{eq:thetas_full_z}
\theta_{z}:=H_{z}-\Pi_{3}H_{z},
\quad
\theta_{z,h}:=H_{z,h}-\Pi_{3}H_{z}.
\end{equation}
\begin{lemma}\label{lemma:identity_dg_1}
There exists a constant $C>0$ independent of $h$ such that
\begin{align*}
\sum_{i=1}^{N_{x}}\Big[
-\int_{I_{i}}[(\hat{\theta}_z\eta_{x,h}^{-})_{x,j+\frac{1}{2}}
-(\hat{\theta}_z\eta_{x,h}^{+})_{x,j-\frac{1}{2}}]dx
+\int_{K_{ij}}\theta_{z}\partial_{y}\eta_{x,h}
\Big]
\le C h^{2k+2}+\|\eta_{x,h}\|^{2},
\\
\sum_{j=1}^{N_{y}}\Big[
\int_{J_{j}}[(\hat{\theta}_z\eta_{y,h}^{-})_{i+\frac{1}{2},y}
-(\hat{\theta}_z\eta_{y,h}^{+})_{i-\frac{1}{2},y}]dy
-\int_{K_{ij}}\theta_{z}\partial_{x}\eta_{y,h}
\Big]
\le C h^{2k+2}+\|\eta_{y,h}\|^{2}.
\end{align*}
\end{lemma}
\proof
See \cite[Lemma~3.4]{Li:17b}.
\close
\begin{lemma}\label{lemma:identity_dg_2}
There exists a constant $C>0$ independent of $h$ such that
\begin{align*}
\sum_{i=1}^{N_{x}}\Big[
-\int_{I_{i}}[(\hat{\theta}_z\eta_{x,h}^{-})_{x,j+\frac{1}{2}}
-(\hat{\theta}_z\eta_{x,h}^{+})_{x,j-\frac{1}{2}}]dx
&+\int_{K_{ij}}\theta_{z}\partial_{y}\eta_{x,h}
\Big]
-\sum_{i=1}^{N_{x}}c_{0}\int_{I_{i}}[\eta_{x,h}^{+}(x,y_{\frac{1}{2}})]^{2}dx\\
&\le C h^{2k+2}+\|\eta_{x,h}\|^{2},
\\
\sum_{j=1}^{N_{y}}\Big[
\int_{J_{j}}[(\hat{\theta}_z\eta_{y,h}^{-})_{i+\frac{1}{2},y}
-(\hat{\theta}_z\eta_{y,h}^{+})_{i-\frac{1}{2},y}]dy
&-\int_{K_{ij}}\theta_{z}\partial_{x}\eta_{y,h}
\Big]
-\sum_{j=1}^{N_{y}}c_{0}\int_{J_{j}}[\eta_{y,h}^{+}(x_{\frac{1}{2}},y)]^{2}dy\\
&\le C h^{2k+2}+\|\eta_{y,h}\|^{2}.
\end{align*}
\end{lemma}
\proof
See \cite[Lemma~3.5]{Li:17b}.
\close
\begin{remark}\label{remark:dg2}
When $c_{0}=0$, we obtain PEC boundary condition without the jump terms in
\eqref{eq:jump6} and \eqref{eq:jump8}.
In this case, we can only control the term
$\sum_{1\le i \le N_{x}}\int_{I_{i}}(\theta_{z}^{+},\eta_{z,h}^{+})(x,c)$ as follows
\begin{align*}
&\sum_{1\le i \le N_{x}}\int_{I_{i}}(\theta_{z}^{+},\eta_{z,h}^{+})(x,c)
\le h^{-1} \int_{r}^{s}(\theta_{z}^{+})^{2}(x,c)
+h \int_{r}^{s}(\eta_{z,h}^{+})^{2}(x,c)
\le C h^{2k+1} + h \|\eta_{z,h}\|^{2},
\end{align*}
by an inverse inequality. Therefore we lose half an order.
\end{remark}
The following result formulates the announced error estimate for the semi-discrete problem.
As in many cases of qualitative estimates, higher regularity requirements are placed
on the weak solution, which of course do not have to be met in all real world situations.
In particular, we assume that the semi-discrete solution is uniformly bounded.
In some special cases, this assumption can be removed at the expense of additional conditions,
mainly a smallness condition to the nonlinearity
\cite[Thm.~4.1]{Bokil:17}, \cite[Thm.~3]{Lyu:21}.
\begin{theorem}\label{th:dg_semi_error_estimate}
Suppose that a weak solution
$(E_{x},E_{y},H_{z})^T \in C^1(0,T,H^{k+1}(\Omega))^3$, $k\in\N$,
of the system \eqref{eq:2dsystem}, and a finite element solution
$(E_{x,h},E_{y,h},H_{z,h})^T\in C^{1}(0,T,U_{h}^{k}\cap L_\infty(\Omega))^3$
of the system \eqref{eq:w_x}--\eqref{eq:w_d2} with the initial data
$E_{x,h}^{0}:=\Pi_{1}E_{x}^{0}$, $E_{y,h}^{0}:=\Pi_{2}E_{y}^{0}$, $H_{z,h}^{0}:=\Pi_{3}H_{z}^{0}$,
respectively exist,
where the $L_\infty$-boundedness of the finite element solution is uniform w.r.t. $h$.
Then, if $h>0$ is sufficiently small, the following error estimate holds
with a coefficient $C(t)>0$ independent of $h$:
\begin{gather*}
\|E_{x}(t)-E_{x,h}(t)\|_{\eps_{0}(1+\chi^{(1)})}+\|E_{y}(t)-E_{y,h}(t)\|_{\eps_{0}(1+\chi^{(1)})}\\
+\|H_{z}(t)-H_{z,h}(t)\|_{\mu_{0}}
\le C(t)h^{k+1},
\quad t\in (0,T).
\end{gather*}
(The concrete structure of $C(t)$ will become apparent from the proof.)
\end{theorem}
\proof
Subtracting the equations \eqref{eq:w_x}--\eqref{eq:w_d2} from the weak formulations
of the equations \eqref{eq:2dsystem}, using the error identities \eqref{eq:zeta_full_x},
\eqref{eq:zeta_full_y}, and \eqref{eq:xi_full_z},
for all test functions $\Phi_{1h}$, $\Phi_{2h}$, $\Phi_{3h}\in Q_{k}(K_{ij})$,
we obtain
\begin{align}
&\int_{K_{ij}}\partial_{t}(D_{x}-D_{x,h})\Phi_{1h}\nonumber\\
&-\int_{I_{i}}[(\hat{\xi_{z}}\Phi_{1h}^{-})_{x,j+\frac{1}{2}}
-(\hat{\xi_{z}}\Phi_{1h}^{+})_{x,j-\frac{1}{2}}]dx\nonumber\\
&+\int_{K_{ij}}\xi_{z}\partial_{y}\Phi_{1h}=0, \label{eq:error_x}\\
&\int_{K_{ij}}\partial_{t}(D_{y}-D_{y,h})\Phi_{2h}
+\int_{J_{j}}[(\hat{\xi_{z}}\Phi_{2h}^{-})_{i+\frac{1}{2},y}\nonumber\\
&-(\hat{\xi_{z}}\Phi_{2h}^{+})_{i-\frac{1}{2},y}]dy\nonumber\\
&-\int_{K_{ij}}\xi_{z}\partial_{x}\Phi_{2h}=0, \label{eq:error_y}\\
&\int_{K_{ij}}\mu_{0}\partial_{t}\xi_{z}\Phi_{3h}
+\int_{J_{j}}[(\hat{\zeta_{y}}\Phi_{3h}^{-})_{i+\frac{1}{2},y}
-(\hat{\zeta_{y}}\Phi_{3h}^{+})_{i-\frac{1}{2},y}]dy\nonumber\\
&-\int_{K_{ij}}\zeta_{y}\partial_{x}\Phi_{3h}
-\int_{I_{i}}[(\hat{\zeta_{x}}\Phi_{3h}^{-})_{x,j+\frac{1}{2}}
-(\hat{\zeta_{x}}\Phi_{3h}^{+})_{x,j-\frac{1}{2}}]dx\nonumber\\
&+\int_{K_{ij}}\zeta_{x}\partial_{y}\Phi_{3h}=0,\nonumber
\end{align}
\begin{align}
&\int_{K_{ij}}\partial_{t}(D_{x}-D_{x,h})\Phi_{1h}
=\int_{K_{ij}}\eps_{0}(1+\chi^{(1)})\partial_{t}\zeta_{x}\Phi_{1h}\nonumber\\
&+\int_{K_{ij}}\eps_{0}\chi^{(3)}\Big[[\vert\EE\vert^2
-\vert\EE_h\vert^2]\partial_{t}E_{x}\Phi_{1h}\nonumber\\
&+\vert\EE_h\vert^2\partial_{t}[E_{x}-E_{x,h}]\Phi_{1h}
+2\big([E_{x}^{2}-E_{x,h}^{2}]\partial_{t}E_{x}\Phi_{1h}\nonumber\\
&+[E_{x}E_{y}-E_{x,h}E_{y,h}]\partial_{t}E_{y}\Phi_{1h}\big)\nonumber\\
&+2\big(E_{x,h}^{2}\partial_{t}[E_{x}-E_{x,h}]\Phi_{1h}\nonumber\\
&+E_{x,h}E_{y,h}\partial_{t}[E_{y}-E_{y,h}]\Phi_{1h}\big)\Big],\label{eq:error_d1}
\end{align}
\begin{align}
&\int_{K_{ij}}\partial_{t}(D_{y}-D_{y,h})\Phi_{2h}
=\int_{K_{ij}}\eps_{0}(1+\chi^{(1)})\partial_{t}\zeta_{y}\Phi_{2h}\nonumber\\
&+\int_{K_{ij}}\eps_{0}\chi^{(3)}\Big[[\vert\EE\vert^2
-\vert\EE_h\vert^2]\partial_{t}E_{y}\Phi_{2h}\nonumber\\
&+\vert\EE_h\vert^2\partial_{t}[E_{y}-E_{y,h}]\Phi_{2h}
+2\big([E_{y}^{2}-E_{y,h}^{2}]\partial_{t}E_{y}\Phi_{2h}\nonumber\\
&+[E_{x}E_{y}-E_{x,h}E_{y,h}]\partial_{t}E_{x}\Phi_{2h}\big)\nonumber\\
&+2\big(E_{y,h}^{2}\partial_{t}[E_{y}-E_{y,h}]\Phi_{2h}\nonumber\\
&+E_{x,h}E_{y,h}\partial_{t}[E_{x}-E_{x,h}]\Phi_{2h}\big)\Big].\label{eq:error_d2}
\end{align}
First we substitute the equations \eqref{eq:error_d1}--\eqref{eq:error_d2} into
the equations \eqref{eq:error_x}--\eqref{eq:error_y}, respectively.
Further decomposing the terms in the resulting equations using
\eqref{eq:etas_full_x}, \eqref{eq:etas_full_y} and \eqref{eq:thetas_full_z} and taking
$\Phi_{1h}:=\eta_{x,h}$, $\Phi_{2h}:=\eta_{y,h}$ and $\Phi_{3h}:=\theta_{z,h}$,
we obtain, after a few slight rearrangements,
\begin{align}
&\int_{K_{ij}}\eps_{0}(1+\chi^{(1)})\partial_t\eta_{x,h}\,\eta_{x,h}
+\int_{K_{ij}}\eps_{0}\chi^{(3)}\Big[\vert\EE_h\vert^2\partial_t\eta_{x,h}\,\eta_{x,h}\nonumber\\
&+2 E_{x,h}^{2}\partial_t\eta_{x,h}\,\eta_{x,h}
+2E_{x,h}E_{y,h}\partial_t\eta_{y,h}\,\eta_{x,h}\Big]\nonumber\\
&-\int_{I_{i}}[(\hat{\theta}_{z,h}\eta_{x,h}^{-})_{x,j+\frac{1}{2}}
-(\hat{\theta}_{z,h}\eta_{x,h}^{+})_{x,j-\frac{1}{2}}]dx
+\int_{K_{ij}}\theta_{z,h}\partial_{y}\eta_{x,h}\nonumber\\
&=\int_{K_{ij}}\eps_{0}(1+\chi^{(1)})\partial_t\eta_{x}\,\eta_{x,h}
+\int_{K_{ij}}\eps_{0}\chi^{(3)}\vert\EE_h\vert^2\partial_t\eta_{x}\,\eta_{x,h}\nonumber\\
&-\int_{I_{i}}[(\hat{\theta}_z\eta_{x,h}^{-})_{x,j+\frac{1}{2}}
-(\hat{\theta}_z\eta_{x,h}^{+})_{x,j-\frac{1}{2}}]dx
+\int_{K_{ij}}\theta_{z}\partial_{y}\eta_{x,h}\nonumber\\
&+\int_{K_{ij}}\eps_{0}\chi^{(3)}\Big[[\vert\EE\vert^2
-\vert\EE_h\vert^2]\partial_t E_{x}\,\eta_{x,h}\nonumber\\
&+2[E_{x}^{2}-E_{x,h}^{2}]\partial_t E_{x}\,\eta_{x,h}+2[E_{x}E_{y}
-E_{x,h}E_{y,h}]\partial_t E_{y}\,\eta_{x,h}\nonumber\\
&+2E_{x,h}^{2}\partial_t\eta_{x}\,\eta_{x,h}+2E_{x,h}E_{y,h}\partial_t\eta_{y}\,\eta_{x,h}\big)\Big],
\label{eq:error_1}
\end{align}
and
\begin{align}
&\int_{K_{ij}}\eps_{0}(1+\chi^{(1)})\partial_t\eta_{y,h}\,\eta_{y,h}
+\int_{K_{ij}}\eps_{0}\chi^{(3)}\Big[\vert\EE_h\vert^2\partial_t\eta_{y,h}\,\eta_{y,h}\nonumber\\
&+2E_{y,h}^{2}\partial_t\eta_{y,h}\,\eta_{y,h}
+2E_{x,h}E_{y,h}\partial_t\eta_{x,h}\,\eta_{y,h}\Big]\nonumber\\
&+\int_{J_{j}}[(\hat{\theta}_{z,h}\eta_{y,h}^{-})_{i+\frac{1}{2},y}
-(\hat{\theta}_{z,h}\eta_{y,h}^{+})_{i-\frac{1}{2},y}]dy
-\int_{K_{ij}}\theta_{z,h}\partial_{x}\eta_{y,h}\nonumber\\
&=\int_{K_{ij}}\eps_{0}(1+\chi^{(1)})\partial_t\eta_{y}\,\eta_{y,h}
+\int_{K_{ij}}\eps_{0}\chi^{(3)}\vert\EE_h\vert^2\partial_t\eta_{y}\,\eta_{y,h}\nonumber\\
&+\int_{J_{j}}[(\hat{\theta}_z\eta_{y,h}^{-})_{i+\frac{1}{2},y}
-(\hat{\theta}_z\eta_{y,h}^{+})_{i-\frac{1}{2},y}]dy
-\int_{K_{ij}}\theta_{z}\partial_{x}\eta_{y,h}\nonumber\\
&+\int_{K_{ij}}\eps_{0}\chi^{(3)}\Big[[\vert\EE\vert^2
-\vert\EE_h\vert^2]\partial_t E_{y}\,\eta_{y,h}\nonumber\\
&+2[E_{y}^{2}-E_{y,h}^{2}]\partial_t E_{y}\,\eta_{y,h}
+2[E_{x}E_{y}-E_{x,h}E_{y,h}]\partial_t E_{x}\,\eta_{y,h}\nonumber\\
&+2E_{y,h}^{2}\partial_t\eta_{y}\,\eta_{y,h}
+2E_{x,h}E_{y,h}\partial_t\eta_{x}\,\eta_{y,h}\Big].
\label{eq:error_2}
\end{align}
For the magnetic field we have that
\begin{align}
&\int_{K_{ij}}\mu_{0}\partial_t\theta_{z,h}\,\theta_{z,h}
+\int_{J_{j}}[(\hat{\eta}_{y,h}\theta_{z,h}^{-})_{i+\frac{1}{2},y}
-(\hat{\eta}_{y,h}\theta_{z,h}^{+})_{i-\frac{1}{2},y}]dy\nonumber\\
&-\int_{I_{i}}[(\hat{\eta}_{x,h}\theta_{z,h}^{-})_{x,j+\frac{1}{2}}
-(\hat{\eta}_{x,h}\theta_{z,h}^{+})_{x,j-\frac{1}{2}}]dx\nonumber\\
&-\int_{K_{ij}}\eta_{y,h}\partial_{x}\theta_{z,h}
+\int_{K_{ij}}\eta_{x,h}\partial_{y}\theta_{z,h}\nonumber\\
&=\int_{K_{ij}}\mu_{0}\partial_t\theta_{z}\,\theta_{z,h}
+\int_{J_{j}}[(\hat{\eta}_y\theta_{z,h}^{-})_{i+\frac{1}{2},y}
-(\hat{\eta}_y\theta_{z,h}^{+})_{i-\frac{1}{2},y}]dy\nonumber\\
&-\int_{K_{ij}}\eta_{y}\partial_{x}\theta_{z,h}
-\int_{I_{i}}[(\hat{\eta}_x\theta_{z,h}^{-})_{x,j+\frac{1}{2}}
-(\hat{\eta}_x\theta_{z,h}^{+})_{x,j-\frac{1}{2}}]dx\nonumber\\
&+\int_{K_{ij}}\eta_{x}\partial_{y}\theta_{z,h}.\label{eq:error_3}
\end{align}
Now we apply similar arguments as in the proof of Thm.~\ref{th:dg_weaKenergy}.
Adding the equations \eqref{eq:error_1}--\eqref{eq:error_3}, summing over
the indices $1\le i\le N_{x}$ and $1\le j\le N_{y}$ and making use of
the identities \eqref{eq:identities12}, we obtain the left-hand side
of the result as
\[ 
LHS = LHSL + LHSN
\] 
with
\begin{equation}\label{eq:error_4l}
\begin{aligned}
LHSL
&:=\frac{1}{2}\frac{d}{dt}\big[\|\eta_{x,h}\|_{\eps_{0}(1+\chi^{(1)})}^{2}
+\|\eta_{y,h}\|_{\eps_{0}(1+\chi^{(1)})}^{2}+\|\theta_{z,h}\|_{\mu_{0}}^{2}\big]\\
&+\sum_{i=1}^{N_{x}}\int_{I_{i}}((\hat{\theta}_{z,h}
-\theta_{z,h}^{+})\eta_{x,h}^{+})(x,y_{\frac{1}{2}})dx\\
&+\sum_{j=1}^{N_{y}}\int_{J_{j}}((\theta_{z,h}^{+}
-\hat{\theta}_{z,h})\eta_{y,h}^{+})(x_{\frac{1}{2}},y)dy\\
&=\frac{1}{2}\frac{d}{dt}\big[\|\eta_{x,h}\|_{\eps_{0}(1+\chi^{(1)})}^{2}
+\|\eta_{y,h}\|_{\eps_{0}(1+\chi^{(1)})}^{2}+\|\theta_{z,h}\|_{\mu_{0}}^{2}\big]\\
&+\sum_{i=1}^{N_{x}} c_{0}\int_{I_{i}}[\eta_{x,h}^{+}(x,y_{\frac{1}{2}})]^{2}dx\\
&+\sum_{j=1}^{N_{y}} c_{0}\int_{J_{j}}[\eta_{y,h}^{+}(x_{\frac{1}{2}},y)]^{2}dy,
\end{aligned}
\end{equation}
where the last equation follows from the definition of the boundary flux densities
\eqref{eq:jump6}, \eqref{eq:jump8} (cf.\ \eqref{eq:identities12}).
Furthermore,
\begin{equation}\label{eq:error_6}
\begin{aligned}
LHSN
&=\int_{\Omega}\eps_{0}\chi^{(3)}\Big[
\frac12 \partial_t\big[\vert\EE_h\vert^2\big(\eta_{x,h}^2+\eta_{y,h}^2\big)\big]\\
&+ \partial_t\big(E_{x,h}\eta_{x,h} + E_{y,h}\eta_{y,h}\big)^2\\
&-\frac12 \big(\eta_{x,h}^2+\eta_{y,h}^2\big)\partial_t\vert\EE_h\vert^2
- \partial_tE_{x,h}^{2}\,\eta_{x,h}^2\\
&- \partial_t E_{y,h}^{2}\,\eta_{y,h}^2
- 2 \partial_t\big(E_{x,h}E_{y,h}\big)\eta_{x,h}\eta_{y,h}\Big].
\end{aligned}
\end{equation}
The right-hand side gets the form
\[ 
RHS = RHSL + RHSN,
\] 
where
\begin{align*}
RHSL
&:=\int_{\Omega}\big[\eps_{0}(1+\chi^{(1)})[\partial_t\eta_{x}\,\eta_{x,h}
+\partial_t\eta_{y}\,\eta_{y,h}]+\mu_{0}\partial_t\theta_{z}\,\theta_{z,h}\big]
\\&
+\sum_{i=1}^{N_{x}}\Big[
-\int_{I_{i}}[(\hat{\theta}_z\eta_{x,h}^{-})_{x,j+\frac{1}{2}}
-(\hat{\theta}_z\eta_{x,h}^{+})_{x,j-\frac{1}{2}}]dx
+\int_{K_{ij}}\theta_{z}\partial_{y}\eta_{x,h}
\Big]\\
&+\sum_{j=1}^{N_{y}}\Big[
\int_{J_{j}}[(\hat{\theta}_z\eta_{y,h}^{-})_{i+\frac{1}{2},y}
-(\hat{\theta}_z\eta_{y,h}^{+})_{i-\frac{1}{2},y}]dy
-\int_{K_{ij}}\theta_{z}\partial_{x}\eta_{y,h}
\Big]
\end{align*}
and
\begin{equation}\label{eq:error_5n}
\begin{aligned}
RHSN
&:= \int_{\Omega}\eps_{0}\chi^{(3)}\Big[
\big[\vert\EE\vert^2-\vert\EE_h\vert^2\big]
\big[\partial_t E_{x}\,\eta_{x,h}+\partial_t E_{y}\,\eta_{y,h}\big]\\
&+2[E_{x}^{2}-E_{x,h}^{2}]\partial_t E_{x}\,\eta_{x,h}
+2[E_{y}^{2}-E_{y,h}^{2}]\partial_t E_{y}\,\eta_{y,h}\\
&+2[E_{x}E_{y}-E_{x,h}E_{y,h}]
\big[\partial_t E_{y}\,\eta_{x,h}+\partial_t E_{x}\,\eta_{y,h}\big]\\
&+2E_{x,h}^{2}\partial_t\eta_{x}\,\eta_{x,h}
+2E_{y,h}^{2}\partial_t\eta_{y}\,\eta_{y,h}\\
&+2E_{x,h}E_{y,h}\big[\partial_t\eta_{y}\,\eta_{x,h}+\partial_t\eta_{x}\,\eta_{y,h}\big]\\
&+\vert\EE_h\vert^2\big[\partial_t\eta_{x}\,\eta_{x,h}+\partial_t\eta_{y}\,\eta_{y,h}\big]
\Big].
\end{aligned}
\end{equation}
Next, using $\partial_{t}E_{x}=\partial_{t}\eta_{x}+\partial_t(\Pi_{1}E_{x})$ and
$\partial_{t}E_{y}=\partial_{t}\eta_{y}+\partial_t(\Pi_{2}E_{y})$
(see \eqref{eq:etas_full_x}, \eqref{eq:etas_full_y})
in the nonlinear terms \eqref{eq:error_5n}, we obtain
\begin{align*}
RHSN
&=\int_{\Omega}\eps_{0}\chi^{(3)}\Big[\big[\big(E_{x}+E_{x,h}\big)\big(E_{x}-E_{x,h}\big)\\
&+\big(E_{y}+E_{y,h}\big)\big(E_{y}-E_{y,h}\big)\big]\big[\partial_t(\Pi_{1} E_{x})\,\eta_{x,h}
+\partial_t(\Pi_{2}E_{y})\,\eta_{y,h}\big]\\
&+\vert\EE\vert^2\big[\partial_t\eta_{x}\,\eta_{x,h}
+\partial_t\eta_{y}\,\eta_{y,h}\big]\\
&+2[\big(E_{x}+E_{x,h}\big)\big(E_{x}-E_{x,h}\big)]\partial_t(\Pi_{1}E_{x})\,\eta_{x,h}\\
&+2E_{x}^{2}\partial_t\eta_{x}\,\eta_{x,h}
+2[\big(E_{y}+E_{y,h}\big)\big(E_{y}-E_{y,h}\big)]\partial_t(\Pi_{2}E_{y})\,\eta_{y,h}\\
&+2E_{y}^{2}\partial_t\eta_{y}\,\eta_{y,h}
+2[E_{y}\big(E_{x}-E_{x,h}\big)\\
&+E_{x,h}\big(E_{y}-E_{y,h}\big)]\big[\partial_t(\Pi_{2}E_{y})\,\eta_{x,h}
+\partial_t(\Pi_{1}E_{x})\,\eta_{y,h}\big]\\
&+2E_{x}E_{y}\big[\partial_t\eta_{y}\,\eta_{x,h}+\partial_t\eta_{x}\,\eta_{y,h}\big]\Big].
\end{align*}
Furthermore, since $E_{x}-E_{x,h}=\eta_{x}-\eta_{x,h}$ and $E_{y}-E_{y,h}=\eta_{y}-\eta_{y,h}$,
we have, after some rearrangement,
\begin{align*}
RHSN
&=\int_{\Omega}\eps_{0}\chi^{(3)}\Big[\big[E_{x}\eta_{x}\partial_t(\Pi_{1} E_{x})
+E_{x,h}\eta_{x}\partial_t(\Pi_{1} E_{x})\\
&+E_{y}\eta_{y}\partial_t(\Pi_{1} E_{x})
+E_{y,h}\eta_{y}\partial_t(\Pi_{1} E_{x})\\
&+\vert\EE\vert^2\partial_t\eta_{x}
+2E_{x}\eta_{x}\partial_t(\Pi_{1}E_{x})+2E_{x,h}\eta_{x}\partial_t(\Pi_{1}E_{x})\\
&+2E_{x}^{2}\partial_t\eta_{x}+2E_{y}\eta_{x}\partial_t(\Pi_{2}E_{y})
+2E_{x,h}\eta_{y}\partial_t(\Pi_{2}E_{y})\\
&+2E_{x}E_{y}\partial_t\eta_{y}\big]\eta_{x,h}\nonumber\\
&+\big[E_{x}\eta_{x}\partial_t(\Pi_{2}E_{y})+E_{x,h}\eta_{x}\partial_t(\Pi_{2}E_{y})\\
&+E_{y}\eta_{y}\partial_t(\Pi_{2}E_{y})+E_{y,h}\eta_{y}\partial_t(\Pi_{2}E_{y})\\
&+\vert\EE\vert^2\partial_t\eta_{y}+2E_{y}\eta_{y}\partial_t(\Pi_{2}E_{y})
+2E_{y,h}\eta_{y}\partial_t(\Pi_{2}E_{y})\\
&+2E_{y}^{2}\partial_t\eta_{y}
+2E_{y}\eta_{x}\partial_t(\Pi_{1}E_{x})+2E_{x,h}\eta_{y}\partial_t(\Pi_{1}E_{x})\\
&+2E_{x}E_{y}\partial_t\eta_{x}\big]\eta_{y,h}\nonumber\\
&+\big[-E_{x}\partial_t(\Pi_{1} E_{x})-E_{x,h}\partial_t(\Pi_{1} E_{x})
-2E_{x}\partial_t(\Pi_{1}E_{x})\\
&-2E_{x,h}\partial_t(\Pi_{1}E_{x})
-2E_{y}\partial_t(\Pi_{2}E_{y})\big]\eta_{x,h}^{2}\nonumber\\
&+\big[-E_{y}\partial_t(\Pi_{2}E_{y})-E_{y,h}\partial_t(\Pi_{2}E_{y})\\
&-2E_{y}\partial_t(\Pi_{2}E_{y})-2E_{y,h}\partial_t(\Pi_{2}E_{y})\\
&-2E_{x,h}\partial_t(\Pi_{1}E_{x})\big]\eta_{y,h}^{2}\nonumber\\
&+\big[-E_{y}\partial_t(\Pi_{1} E_{x})-E_{y,h}\partial_t(\Pi_{1} E_{x})\\
&-E_{x}\partial_t(\Pi_{2}E_{y})-E_{x,h}\partial_t(\Pi_{2}E_{y})\\
&-2E_{y}\partial_t(\Pi_{1}E_{x})
-2E_{x,h}\partial_t(\Pi_{2}E_{y})\big]\eta_{x,h}\eta_{y,h}\Big].
\end{align*}
In a next step, we shift
the last two terms of \eqref{eq:error_4l} to $RHSL$ and
the last four terms of \eqref{eq:error_6} to $RHSN$.
Then the new left-hand side is
\begin{equation}\label{eq:error_4lhsnew}
\begin{aligned}
LHS'
&:=\frac{1}{2}\frac{d}{dt}\big[\|\eta_{x,h}\|_{\eps_{0}(1+\chi^{(1)})}^{2}
+\|\eta_{y,h}\|_{\eps_{0}(1+\chi^{(1)})}^{2}+\|\theta_{z,h}\|_{\mu_{0}}^{2}\big]\\
&+\int_{\Omega}\eps_{0}\chi^{(3)}\Big[
\frac12 \partial_t\big[\vert\EE_h\vert^2\big(\eta_{x,h}^2+\eta_{y,h}^2\big)\big]\\
&+ \partial_t\big(E_{x,h}\eta_{x,h} + E_{y,h}\eta_{y,h}\big)^2\Big],
\end{aligned}
\end{equation}
whereas the new right-hand side is
\[
RHS':=RHSL'+RHSN'
\]
with
\begin{align*}
RHSL':=&RHSL-\sum_{i=1}^{N_{x}} c_{0}\int_{I_{i}}[\eta_{x,h}^{+}(x,y_{\frac{1}{2}})]^{2}dx
-\sum_{j=1}^{N_{y}} c_{0}\int_{J_{j}}[\eta_{y,h}^{+}(x_{\frac{1}{2}},y)]^{2}dy,\\
RHSN':=&RHSN+\int_{\Omega}\eps_{0}\chi^{(3)}\Big[
\frac12 \big(\eta_{x,h}^2+\eta_{y,h}^2\big)\partial_t\vert\EE_h\vert^2\\
&\qquad\qquad +\partial_tE_{x,h}^{2}\,\eta_{x,h}^2
+\partial_t E_{y,h}^{2}\,\eta_{y,h}^2
+ 2 \partial_t\big(E_{x,h}E_{y,h}\big)\eta_{x,h}\eta_{y,h}\Big].
\end{align*}
The first three terms from $RHSL'$
are estimated using the Cauchy-Schwarz inequality and Lemma \ref{lemma:pro_error_dg}:
\begin{align*}
&\int_{\Omega}\big[\eps_{0}(1+\chi^{(1)})[\partial_t\eta_{x}\,\eta_{x,h}
+\partial_t\eta_{y}\,\eta_{y,h}]+\mu_{0}\partial_t\theta_{z}\,\theta_{z,h}\big]\\
&\le \|\partial_{t}\eta_{x}\|_{\eps_{0}(1+\chi^{(1)})}\|\eta_{x,h}\|_{\eps_{0}(1+\chi^{(1)})}\\
&+\|\partial_{t}\eta_{y}\|_{\eps_{0}(1+\chi^{(1)})}\|\eta_{y,h}\|_{\eps_{0}(1+\chi^{(1)})}
+\|\partial_{t}\theta_{z}\|_{\mu_{0}}\|\theta_{z,h}\|_{\mu_{0}}\\
&\le C_{1}h^{k+1}[\|\eta_{x,h}\|_{\eps_{0}(1+\chi^{(1)})}+\|\eta_{y,h}\|_{\eps_{0}(1+\chi^{(1)})}
+\|\theta_{z,h}\|_{\mu_{0}}],
\end{align*}
where the constant $C_{1}>0$ depends on
$\|\eps_{0}(1+\chi^{(1)})\|_{L_\infty(\Omega)}$,
$\|\mu_{0}\|_{L_\infty(\Omega)}$,
$\|\partial_{t}E_{x}\|_{H^{k+1}(\Omega)}$,
$\|\partial_{t}E_{y}\|_{H^{k+1}(\Omega)}$,
and
$\|\partial_{t}H_{z}\|_{H^{k+1}(\Omega)}$,
as can be seen be the following exemplary argument:
\begin{align*}
\|\partial_{t}\eta_{x}\|_{\eps_{0}(1+\chi^{(1)})}
\le \|\eps_{0}(1+\chi^{(1)})\|_{L_\infty(\Omega)}\|\partial_{t}\eta_{x}\|
\le C_{1}h^{k+1}\|\partial_{t}E_{x}\|_{H^{k+1}(\Omega)}.
\end{align*}
The remaining terms from $RHSL'$ are
estimated by means of the Lemmata \ref{lemma:identity_dg_1}, \ref{lemma:identity_dg_2}:
\begin{align*}
&\sum_{i=1}^{N_{x}}\Big[
-\int_{I_{i}}[(\hat{\theta}_z\eta_{x,h}^{-})_{x,j+\frac{1}{2}}
-(\hat{\theta}_z\eta_{x,h}^{+})_{x,j-\frac{1}{2}}]dx
+\int_{K_{ij}}\theta_{z}\partial_{y}\eta_{x,h}
\Big]\\
&+\sum_{j=1}^{N_{y}}\Big[
\int_{J_{j}}[(\hat{\theta}_z\eta_{y,h}^{-})_{i+\frac{1}{2},y}
-(\hat{\theta}_z\eta_{y,h}^{+})_{i-\frac{1}{2},y}]dy
-\int_{K_{ij}}\theta_{z}\partial_{x}\eta_{y,h}
\Big]\\
&-\sum_{i=1}^{N_{x}} c_{0}\int_{I_{i}}[\eta_{x,h}^{+}(x,y_{\frac{1}{2}})]^{2}dx
-\sum_{j=1}^{N_{y}} c_{0}\int_{J_{j}}[\eta_{y,h}^{+}(x_{\frac{1}{2}},y)]^{2}dy\\
&\le C h^{2k+2} + \|\eta_{x,h}\|^{2} + \|\eta_{y,h}\|^{2}\\
&\le C h^{2k+2}
+ \|(\eps_{0}(1+\chi^{(1)}))^{-1}\|_{L_\infty(\Omega)}
\big[\|\eta_{x,h}\|_{\eps_{0}(1+\chi^{(1)})}^2
+ \|\eta_{y,h}\|_{\eps_{0}(1+\chi^{(1)})}^2\big].
\end{align*}
The terms from the right-hand side part $RHSN$ can be estimated as follows:
\begin{align*}
RHSN
&\le \|\chi^{(3)}(1+\chi^{(1)})^{-1}\|_{L_\infty(\Omega)}\\
&\times\bigg[
\Big[\|E_{x}\|_{L_\infty(\Omega)}
\|\partial_t(\Pi_{1}E_{x})\|_{L_\infty(\Omega)}\|\eta_{x}\|_{\eps_{0}(1+\chi^{(1)})}\\
&+\|E_{x,h}\|_{L_\infty(\Omega)}\|\eta_{x}\|_{\eps_{0}(1+\chi^{(1)})}
\|\partial_t(\Pi_{1}E_{x})\|_{L_\infty(\Omega)}\\
&+\|E_{y}\|_{L_\infty(\Omega)}\|\eta_{y}\|_{\eps_{0}(1+\chi^{(1)})}
\|\partial_t(\Pi_{1}E_{x})\|_{L_\infty(\Omega)}\\
&+\|E_{y,h}\|_{L_\infty(\Omega)}\|\eta_{y}\|_{\eps_{0}(1+\chi^{(1)})}
\|\partial_t(\Pi_{1}E_{x})\|_{L_\infty(\Omega)}\\
&+\|\vert\EE\vert^2\|_{L_\infty(\Omega)}\|\partial_{t}\eta_{x}\|_{\eps_{0}(1+\chi^{(1)})}\\
&+2\|E_{x}\|_{L_\infty(\Omega)}\|\eta_{x}\|_{\eps_{0}(1+\chi^{(1)})}
\|\partial_t(\Pi_{1}E_{x})\|_{L_\infty(\Omega)}\\
&+2\|E_{x,h}\|_{L_\infty(\Omega)}\|\eta_{x}\|_{\eps_{0}(1+\chi^{(1)})}
\|\partial_t(\Pi_{1}E_{x})\|_{L_\infty(\Omega)}\\
&+2\|E_{x}^{2}\|_{L_\infty(\Omega)}\|\partial_{t}\eta_{x}\|_{\eps_{0}(1+\chi^{(1)})}\\
&+2\|E_{y}\|_{L_\infty(\Omega)}\|\eta_{x}\|_{\eps_{0}(1+\chi^{(1)})}
\|\partial_t(\Pi_{2}E_{y})\|_{L_\infty(\Omega)}\\
&+2\|E_{x,h}\|_{L_\infty(\Omega)}\|\eta_{y}\|_{\eps_{0}(1+\chi^{(1)})}
\|\partial_t(\Pi_{2}E_{y})\|_{L_\infty(\Omega)}\\
&+2\|E_{x}\|_{L_\infty(\Omega)}\|E_{y}\|_{L_\infty(\Omega)}
\|\partial_{t}\eta_{y}\|_{\eps_{0}(1+\chi^{(1)})}\Big]
\|\eta_{x,h}\|_{\eps_{0}(1+\chi^{(1)})}\\
&+\Big[\|E_{x}\|_{L_\infty(\Omega)}\|\eta_{x}\|_{\eps_{0}(1+\chi^{(1)})}
\|\partial_t(\Pi_{2}E_{y})\|_{L_\infty(\Omega)}\\
&+\|E_{x,h}\|_{L_\infty(\Omega)}\|\eta_{x}\|_{\eps_{0}(1+\chi^{(1)})}
\|\partial_t(\Pi_{2}E_{y})\|_{L_\infty(\Omega)}\\
&+\|E_{y}\|_{L_\infty(\Omega)}\|\eta_{y}\|_{\eps_{0}(1+\chi^{(1)})}
\|\partial_t(\Pi_{2}E_{y})\|_{L_\infty(\Omega)}\\
&+\|E_{y,h}\|_{L_\infty(\Omega)}\|\eta_{y}\|_{\eps_{0}(1+\chi^{(1)})}
\|\partial_t(\Pi_{2}E_{y})\|_{L_\infty(\Omega)}\\
&+\|\vert\EE\vert^2\|_{L_\infty(\Omega)}\|\partial_{t}\eta_{y}\|_{\eps_{0}(1+\chi^{(1)})}\\
&+2\|E_{y}\|_{L_\infty(\Omega)}\|\eta_{y}\|_{\eps_{0}(1+\chi^{(1)})}
\|\partial_t(\Pi_{2}E_{y})\|_{L_\infty(\Omega)}\\
&+2\|E_{y,h}\|_{L_\infty(\Omega)}\|\eta_{y}\|_{\eps_{0}(1+\chi^{(1)})}
\|\partial_t(\Pi_{2}E_{y})\|_{L_\infty(\Omega)}\\
&+2\|E_{y}^{2}\|_{L_\infty(\Omega)}\|\partial_{t}\eta_{y}\|_{\eps_{0}(1+\chi^{(1)})}\\
&+\|E_{y}\|_{L_\infty(\Omega)}\|\eta_{x}\|_{\eps_{0}(1+\chi^{(1)})}
\|\partial_t(\Pi_{1}E_{x})\|_{L_\infty(\Omega)}\\
&+2\|E_{x,h}\|_{L_\infty(\Omega)}\|\eta_{y}\|_{\eps_{0}(1+\chi^{(1)})}
\|\partial_t(\Pi_{1}E_{x})\|_{L_\infty(\Omega)}\\
&+2\|E_{x}\|_{L_\infty(\Omega)}\|E_{y}\|_{L_\infty(\Omega)}
\|\partial_{t}\eta_{x}\|_{\eps_{0}(1+\chi^{(1)})}\Big]
\|\eta_{y,h}\|_{\eps_{0}(1+\chi^{(1)})}\\
&+\Big[\|E_{x}\|_{L_\infty(\Omega)}\|\partial_t(\Pi_{1}E_{x})\|_{L_\infty(\Omega)}\\
&+\|E_{x,h}\|_{L_\infty(\Omega)}\|\partial_t(\Pi_{1}E_{x})\|_{L_\infty(\Omega)}\\
&+2\|E_{x}\|_{L_\infty(\Omega)}\|\partial_t(\Pi_{1}E_{x})\|_{L_\infty(\Omega)}\\
&+2\|E_{x,h}\|_{L_\infty(\Omega)}\|\partial_t(\Pi_{1}E_{x})\|_{L_\infty(\Omega)}\\
&+2\|E_{y}\|_{L_\infty(\Omega)}\|\partial_t(\Pi_{2}E_{y})\|_{L_\infty(\Omega)}\Big]
\|\eta_{x,h}\|_{\eps_{0}(1+\chi^{(1)})}^2\\
&+\Big[\|E_{y}\|_{L_\infty(\Omega)}\|\partial_t(\Pi_{2}E_{y})\|_{L_\infty(\Omega)}\\
&+\|E_{y,h}\|_{L_\infty(\Omega)}\|\partial_t(\Pi_{2}E_{y})\|_{L_\infty(\Omega)}\\
&+2\|E_{y}\|_{L_\infty(\Omega)}\|\partial_t(\Pi_{2}E_{y})\|_{L_\infty(\Omega)}\\
&+2\|E_{y,h}\|_{L_\infty(\Omega)}\|\partial_t(\Pi_{2}E_{y})\|_{L_\infty(\Omega)}\\
&+2\|E_{x,h}\|_{L_\infty(\Omega)}\|\partial_t(\Pi_{1}E_{x})\|_{L_\infty(\Omega)}\Big]
\|\eta_{y,h}\|_{\eps_{0}(1+\chi^{(1)})}^2\\
&+\Big[\|E_{y}\|_{L_\infty(\Omega)}\|\partial_t(\Pi_{1}E_{x})\|_{L_\infty(\Omega)}\\
&+\|E_{y,h}\|_{L_\infty(\Omega)}\|\partial_t(\Pi_{1}E_{x})\|_{L_\infty(\Omega)}\\
&+\|E_{x}\|_{L_\infty(\Omega)}\|\partial_t(\Pi_{2}E_{y})\|_{L_\infty(\Omega)}\\
&+\|E_{x,h}\|_{L_\infty(\Omega)}\|\partial_t(\Pi_{2}E_{y})\|_{L_\infty(\Omega)}\\
&+2\|E_{y}\|_{L_\infty(\Omega)}\|\partial_t(\Pi_{1}E_{x})\|_{L_\infty(\Omega)}\\
&+2\|E_{x,h}\|_{L_\infty(\Omega)}\|\partial_t(\Pi_{2}E_{y})\|_{L_\infty(\Omega)}\Big]\\
&\times\|\eta_{x,h}\|_{\eps_{0}(1+\chi^{(1)})}\|\eta_{y,h}\|_{\eps_{0}(1+\chi^{(1)})}
\bigg].
\end{align*}
This estimate shows that we have to discuss upper bounds for the terms
$\|E_{x}\|_{L_\infty(\Omega)}$,
$\|E_{y}\|_{L_\infty(\Omega)}$,
$\|\vert\EE\vert^2\|_{L_\infty(\Omega)}$,
$\|E_{x}^{2}\|_{L_\infty(\Omega)}$,
$\|E_{y}^{2}\|_{L_\infty(\Omega)}$,
$\|\partial_t(\Pi_{1}E_{x})\|_{L_\infty(\Omega)}$,
$\|\partial_t(\Pi_{2}E_{y})\|_{L_\infty(\Omega)}$,
$\|\eta_{x}\|_{\eps_{0}(1+\chi^{(1)})}$,
$\|\eta_{y}\|_{\eps_{0}(1+\chi^{(1)})}$,
\linebreak
$\|\partial_{t}\eta_{x}\|_{\eps_{0}(1+\chi^{(1)})}$,
$\|\partial_{t}\eta_{y}\|_{\eps_{0}(1+\chi^{(1)})}$,
$\|E_{x,h}\|_{L_\infty(\Omega)}$,
and
$\|E_{y,h}\|_{L_\infty(\Omega)}.$
The first five terms are bounded thanks to the assumption w.r.t.\ the weak solution
and the continuous embedding $H^{k+1}(\Omega)\subset L_\infty(\Omega)$ for $k\in\N$,
see, e.g., \cite[(3.1.4)]{Ciarlet:02b} (this embedding remains valid for $d=3$, too).
The eighth to eleventh terms are estimated by means of Lemma~\ref{lemma:pro_error_dg},
for instance:
\begin{align*}
\|\eta_{x}\|_{\eps_{0}(1+\chi^{(1)})}
\le \|\eps_{0}(1+\chi^{(1)})\|_{L_\infty(\Omega)}^{1/2}\|\eta_{x}\|
\le C h^{k+1}\|E_x\|_{H^{k+1}(\Omega)},
\end{align*}
where here the constant $C>0$ depends on $\|\eps_{0}(1+\chi^{(1)})\|_{L_\infty(\Omega)}$.
The last two terms are bounded thanks to the assumption w.r.t.\ the numerical solution.

So it remains to investigate the sixth and seventh terms. Taking into account
the commutation property $\partial_t(\Pi_{1}E_{x})=\Pi_{1}(\partial_tE_{x})$,
we first observe that there exist at least one element $K_{ij}$ such that
\[
\|\partial_t(\Pi_{1}E_{x})\|_{L_\infty(\Omega)}
= \|\Pi_{1}(\partial_tE_{x})\|_{L_\infty(\Omega)}
= \|\Pi_{1}(\partial_tE_{x})\|_{L_\infty(K_{ij})}
\]
The latter norm can be estimated by an inverse inequality \cite[Thm.~3.2.6]{Ciarlet:02b}:
\[
\|\Pi_{1}(\partial_tE_{x})\|_{L_\infty(K_{ij})}
\le \vert K_{ij}\vert^{-1/2}\|\Pi_{1}(\partial_tE_{x})\|_{L_2(K_{ij})}
\]
(note that we need only a local variant,
i.e.\ we may omit the inverse assumption \cite[(3.2.28)]{Ciarlet:02b}).
Using the triangle inequality, we get
\begin{align*}
\|\partial_t(\Pi_{1}E_{x})\|_{L_\infty(\Omega)}
&\le \vert K_{ij}\vert^{-1/2}\big[
\|\partial_tE_{x}\|_{L_2(K_{ij})}
+\|\Pi_{1}(\partial_tE_{x})-\partial_tE_{x}\|_{L_2(K_{ij})}
\big]\\
&\le \vert K_{ij}\vert^{-1/2}\big[
\vert K_{ij}\vert^{1/2}\|\partial_tE_{x}\|_{L_\infty(K_{ij})}
+ C \vert K_{ij}\vert^{(k+1)/2}\|\partial_tE_x\|_{H^{k+1}(K_{ij})}
\big].
\end{align*}
The estimate of the first term in the square brackets results from H\"{o}lder's
inequality, whereas the second term is estimated by means of a local variant
of Lemma~\ref{lemma:pro_error_dg}, see \cite[Lemma~3.2]{Cockburn:01c}.
So if the mesh size $h$ is sufficiently small, we get
\begin{align*}
\|\Pi_{1}(\partial_tE_{x})\|_{L_\infty(\Omega)}
&\le \|\partial_tE_{x}\|_{L_\infty(K_{ij})}
+ C \vert K_{ij}\vert^{k/2}\|\partial_tE_x\|_{H^{k+1}(K_{ij})}\\
&\le C \|\partial_tE_x\|_{H^{k+1}(\Omega)},
\end{align*}
where we have used the continuous embedding $H^{k+1}(\Omega)\subset L_\infty(\Omega)$
in the the last step again.
An analogous argument applies to the seventh term.

So in summary we arrive at the estimate
\begin{align*}
RHSN & \le C_{2}h^{k+1}\big[\|\eta_{x,h}\|_{\eps_{0}(1+\chi^{(1)})}
+\|\eta_{y,h}\|_{\eps_{0}(1+\chi^{(1)})}\big]\\
& + C_{3}\|\eta_{x,h}\|_{\eps_{0}(1+\chi^{(1)})}^2
+ C_{4}\|\eta_{y,h}\|_{\eps_{0}(1+\chi^{(1)})}^2\\
&+ C_{5}\|\eta_{x,h}\|_{\eps_{0}(1+\chi^{(1)})}\|\eta_{y,h}\|_{\eps_{0}(1+\chi^{(1)})}\\
& \le C_{2}h^{k+1}\big[\|\eta_{x,h}\|_{\eps_{0}(1+\chi^{(1)})}
+\|\eta_{y,h}\|_{\eps_{0}(1+\chi^{(1)})}\big]\\
&+ C_{6}\big[\|\eta_{x,h}\|_{\eps_{0}(1+\chi^{(1)})}^2
+\|\eta_{y,h}\|_{\eps_{0}(1+\chi^{(1)})}^2\big],
\end{align*}
where the constants $C_{2},C_{6}$ depend on
$\|\chi^{(3)}(1+\chi^{(1)})^{-1}\|_{L_\infty(\Omega)}$,
the $C^1(0,T,H^{k+1}(\Omega))$-norms of $E_x,E_y$
and the $C^{1}(0,T,U_{h}^{k}\cap L_\infty(\Omega))$-norms of $E_{x,h},E_{y,h}$.

Furthermore the remaining terms from $RHSN'$ can be bounded from above by
\begin{align*}
\|\chi^{(3)}(1+\chi^{(1)})^{-1}\|_{L_\infty(\Omega)}
&
\bigg[
\frac{1}{2}\|\partial_{t}\vert\EE_h\vert^2\|_{L_\infty(\Omega)}
\big[\|\eta_{x,h}\|_{\eps_{0}(1+\chi^{(1)})}^{2}+\|\eta_{y,h}\|_{\eps_{0}(1+\chi^{(1)})}^{2}\big]\\
&+\|\partial_{t}E_{x,h}^{2}\|_{L_\infty(\Omega)}\|\eta_{x,h}\|_{\eps_{0}(1+\chi^{(1)})}^{2}
+\|\partial_{t}E_{y,h}^{2}\|_{L_\infty(\Omega)}\|\eta_{y,h}\|_{\eps_{0}(1+\chi^{(1)})}^{2}\\
&+2\|\partial_{t}(E_{x,h}E_{y,h})\|_{L_\infty(\Omega)}
\|\eta_{x,h}\|_{\eps_{0}(1+\chi^{(1)})}\|\eta_{y,h}\|_{\eps_{0}(1+\chi^{(1)})}\bigg].
\end{align*}
Here we have to take care of
$\|\partial_{t}\vert\EE_h\vert^2\|_{L_\infty(\Omega)}$,
$\|\partial_{t}E_{x,h}^{2}\|_{L_\infty(\Omega)}$,
$\|\partial_{t}E_{y,h}^{2}\|_{L_\infty(\Omega)}$,
and
$\|\partial_{t}(E_{x,h}E_{y,h})\|_{L_\infty(\Omega)}$,
but all these terms can be bounded from above by the
$C^{1}(0,T,U_{h}^{k}\cap L_\infty(\Omega))$-norms of $E_{x,h},E_{y,h}$.
Therefore we get the upper bound
\[
C\big[\|\eta_{x,h}\|_{\eps_{0}(1+\chi^{(1)})}^2
+\|\eta_{y,h}\|_{\eps_{0}(1+\chi^{(1)})}^2\big],
\]
where the constant $C$ depends on
$\|\chi^{(3)}(1+\chi^{(1)})^{-1}\|_{L_\infty(\Omega)}$
and the $C^{1}(0,T,U_{h}^{k}\cap L_\infty(\Omega))$-norms of $E_{x,h},E_{y,h}$.
Since such a term already occurs in the upper bound of $RHSN$, we modify
the constant $C_{6}$ correspondingly and conclude
\begin{equation}\label{eq:error_13}
\begin{aligned}
RHSN'
& \le C_{2}h^{k+1}\big[\|\eta_{x,h}\|_{\eps_{0}(1+\chi^{(1)})}
+\|\eta_{y,h}\|_{\eps_{0}(1+\chi^{(1)})}\big]\\
&+ C_{6}\big[\|\eta_{x,h}\|_{\eps_{0}(1+\chi^{(1)})}^2
+\|\eta_{y,h}\|_{\eps_{0}(1+\chi^{(1)})}^2\big].
\end{aligned}
\end{equation}
Combining the right-hand side estimate \eqref{eq:error_13}
with the left the-hand side \eqref{eq:error_4lhsnew}, we obtain
\begin{align*}
&\frac{1}{2}\frac{d}{dt}\big[\|\eta_{x,h}\|_{\eps_{0}(1+\chi^{(1)})}^{2}
+\|\eta_{y,h}\|_{\eps_{0}(1+\chi^{(1)})}^{2}+\|\theta_{z,h}\|_{\mu_{0}}^{2}\big]\\
&
+\int_{\Omega}\eps_{0}\chi^{(3)}\Big[
\frac12 \partial_t\big[\vert\EE_h\vert^2\big(\eta_{x,h}^2+\eta_{y,h}^2\big)\big]
+ \partial_t\big(E_{x,h}\eta_{x,h} + E_{y,h}\eta_{y,h}\big)^2\Big]\\
& \le C_{2}h^{k+1}\big[\|\eta_{x,h}\|_{\eps_{0}(1+\chi^{(1)})}
+\|\eta_{y,h}\|_{\eps_{0}(1+\chi^{(1)})}\big]\\
&+ C_{6}\big[\|\eta_{x,h}\|_{\eps_{0}(1+\chi^{(1)})}^2
+\|\eta_{y,h}\|_{\eps_{0}(1+\chi^{(1)})}^2\big].
\end{align*}
Setting
\begin{align*}
\mathcal{D}_h^2(t)
&:=\|\eta_{x,h}(t)\|_{\eps_{0}(1+\chi^{(1)})}^{2}
+\|\eta_{y,h}(t)\|_{\eps_{0}(1+\chi^{(1)})}^{2}+\|\theta_{z,h}(t)\|_{\mu_{0}}^{2}\\
&+\int_{\Omega}\eps_{0}\chi^{(3)}\Big[
\vert\EE_h(t)\vert^2\big(\eta_{x,h}^2(t)+\eta_{y,h}^2(t)\big)\\
&+ 2\big(E_{x,h}(t)\eta_{x,h}(t) + E_{y,h}(t)\eta_{y,h}(t)\big)^2\Big],
\end{align*}
we get
\begin{align*}
\frac{1}{2}\frac{d}{dt}\mathcal{D}_h^2(t)
& \le C_{2}h^{k+1}\big[\|\eta_{x,h}\|_{\eps_{0}(1+\chi^{(1)})}
+\|\eta_{y,h}\|_{\eps_{0}(1+\chi^{(1)})}\big]\\
&+ C_{6}\big[\|\eta_{x,h}\|_{\eps_{0}(1+\chi^{(1)})}^2
+\|\eta_{y,h}\|_{\eps_{0}(1+\chi^{(1)})}^2\big]\\
& \le C_{2}\sqrt{2}h^{k+1}\mathcal{D}_h(t)
+ C_{6}\mathcal{D}_h^2(t).
\end{align*}
Integrating this inequality with respect to time, we obtain
\[
\mathcal{D}_h^2(t) \le \mathcal{D}_h(0)
+ 2\int_0^t\big[C_{2}\sqrt{2}h^{k+1}\mathcal{D}_h(s) + C_{6}\mathcal{D}_h^2(s)\big]ds.
\]
Now we apply a Gronwall-type lemma \cite[Lemma 4.1]{Dafermos:79} and obtain
\[
\mathcal{D}_h(t) \le \mathcal{D}_h(0)e^{C_{6}t} + C_{2}\sqrt{2}h^{k+1} t e^{C_{6}t}.
\]
From this and the triangle inequality in conjunction with Lemma~\ref{lemma:pro_error_dg}
the statement follows.
\close

\bigskip
\section{The Fully Discrete Scheme}
We divide the time interval $(0,T)$ into $N\in\N$ equally spaced subintervals
by using the nodal points
$t^{n}:=n\Delta t$, $n=0,1,2,\ldots,N$, and $\Delta t:=\frac{T}{N}$.
Given initial values
$(E_{x,h}^{0},E_{y,h}^{0},H_{z,h}^{0})^T\in (U_{h}^{k})^3$
of the electric and magnetic field intensities,
the fully discrete scheme w.r.t.\ the electric and magnetic field intensities
$(E_{x,h}^{n+1},E_{y,h}^{n+1},H_{z,h}^{n+\frac{3}{2}})^T\in (U_{h}^{k})^3$, $n=1,2,\ldots,N-1$,
reads as
\begin{align}
&\int_{K_{ij}}\frac{D_{x,h}^{n+1}-D_{x,h}^{n}}{\Delta t} \Phi_{1h}\nonumber\\
&-\int_{I_{i}}[(\hat{H}_{z,h}^{n+\frac{1}{2}}\Phi_{1h}^{-})_{x,j+\frac{1}{2}}
-(\hat{H}_{z,h}^{n+\frac{1}{2}}\Phi_{1h}^{+})_{x,j-\frac{1}{2}}]dx\nonumber\\
&+\int_{K_{ij}}H_{z,h}^{n+\frac{1}{2}}\partial_{y}\Phi_{1h}
-\int_{K_{ij}}J_{x,h}^{n+\frac{1}{2}}\Phi_{1h}=0, \label{eq:f_x}\\
&\int_{K_{ij}}\frac{D_{y,h}^{n+1}-D_{y,h}^{n}}{\Delta t}\nonumber\\
&+\int_{J_{j}}[(\hat{H}_{z,h}^{n+\frac{1}{2}}\Phi_{2h}^{-})_{i+\frac{1}{2},y}
-(\hat{H}_{z,h}^{n+\frac{1}{2}}\Phi_{2h}^{+})_{i-\frac{1}{2},y}]dy\nonumber\\
&-\int_{K_{ij}}H_{z,h}^{n+\frac{1}{2}}\partial_{x}\Phi_{2h}
-\int_{K_{ij}}J_{y,h}^{n+\frac{1}{2}}\Phi_{2h}=0, \label{eq:f_y}\\
&\int_{K_{ij}}\mu_{0}\frac{H_{z,h}^{n+\frac{3}{2}}
-H_{z,h}^{n+\frac{1}{2}}}{\Delta t}\Phi_{3h}\nonumber\\
&+\int_{J_{j}}[(\hat{E}_{y,h}^{n+1}\Phi_{3h}^{-})_{i+\frac{1}{2},y}
-(\hat{E}_{y,h}^{n+1}\Phi_{3h}^{+})_{i-\frac{1}{2},y}]dy\nonumber\\
&-\int_{K_{ij}}E_{y,h}^{n+1}\partial_{x}\Phi_{3h}
-\int_{I_{i}}[(\hat{E}_{x,h}^{n+1}\Phi_{3h}^{-})_{x,j+\frac{1}{2}}\nonumber\\
&-(\hat{E}_{x,h}^{n+1}\Phi_{3h}^{+})_{x,j-\frac{1}{2}}]dx
+\int_{K_{ij}}E_{x,h}^{n+1}\partial_{y}\Phi_{3h}=0,\label{eq:f_h}\\
&\int_{K_{ij}}(D_{x,h}^{n+1}-D_{x,h}^{n})\Phi_{1h}\nonumber\\
&=\int_{K_{ij}}\eps_{0}(1+\chi^{(1)})(E_{x,h}^{n+1}-E_{x,h}^{n})\Phi_{1h}\nonumber\\
&+\int_{K_{ij}}\eps_{0}\chi^{(3)}\Big[\frac{1}{2}\big[(E_{x,h}^{n+1})^{2}
+(E_{x,h}^{n})^{2}+(E_{y,h}^{n+1})^{2}\nonumber\\
&+(E_{y,h}^{n})^{2}\big](E_{x,h}^{n+1}-E_{x,h}^{n})\Phi_{1h}\label{eq:f_d1}\\
&+\big([(E_{x,h}^{n+1})^{2}+(E_{x,h}^{n})^{2}](E_{x,h}^{n+1}-E_{x,h}^{n})\Phi_{1h}\nonumber\\
&+[E_{x,h}^{n+1}E_{y,h}^{n+1}
+E_{x,h}^{n}E_{y,h}^{n}](E_{y,h}^{n+1}-E_{y,h}^{n})\Phi_{1h}\big)\Big],\nonumber\\
&\int_{K_{ij}}(D_{y,h}^{n+1}-D_{y,h}^{n})\Phi_{2h}\nonumber\\
&=\int_{K_{ij}}\eps_{0}(1+\chi^{(1)})(E_{y,h}^{n+1}-E_{y,h}^{n})\Phi_{2h}\nonumber\\
&+\int_{K_{ij}}\eps_{0}\chi^{(3)}\Big[\frac{1}{2}\big((E_{x,h}^{n+1})^{2}
+(E_{x,h}^{n})^{2}+(E_{y,h}^{n+1})^{2}\nonumber\\
&+(E_{y,h}^{n})^{2}\big)(E_{y,h}^{n+1}-E_{y,h}^{n})\Phi_{2h}
\label{eq:f_d2}\\
&+\big([(E_{y,h}^{n+1})^{2}+(E_{y,h}^{n})^{2}](E_{y,h}^{n+1}-E_{y,h}^{n})\Phi_{2h}\nonumber\\
&+[E_{x,h}^{n+1}E_{y,h}^{n+1}+E_{x,h}^{n}E_{y,h}^{n}](E_{x,h}^{n+1}
-E_{x,h}^{n})\Phi_{2h}\big)\Big]\nonumber
\end{align}
for all test functions $(\Phi_{1h},\Phi_{2h},\Phi_{3h})^T\in (U_{h}^{k})^3$.
The differences $D_{x,h}^{n+1}-D_{x,h}^{n}$ and $D_{y,h}^{n+1}-D_{y,h}^{n}$
play the role of auxiliary variables, and the flux densities are defined by
\begin{align}
\hat{E}_{x,h}^{n+1}(x,y_{j+\frac{1}{2}})
&:=E_{x,h}^{n+1, +}(x,y_{j+\frac{1}{2}})\nonumber\\
&\quad\text{for all } j = 1, 2, 3,\ldots ,N_{y}-1,\label{eq:f_jump1}\\
\hat{E}_{x,h}^{n+1}(x,y_{\frac{1}{2}})
&:=\hat{E}_{x,h}^{n+1}(x,y_{N_{y}+\frac{1}{2}}):=0,\label{eq:f_jump2}\\
\hat{E}_{y,h}^{n+1}(x_{i+\frac{1}{2}},y)
&:=E_{y,h}^{n+1,+}(x_{i+\frac{1}{2}},y)\nonumber\\
&\quad\text{for all } i= 1, 2, 3,\ldots, N_{x}-1,\label{eq:f_jump3}\\
\hat{E}_{y,h}^{n+1}(x_{\frac{1}{2}},y)
&:=\hat{E}_{y,h}^{n+1}(x_{N_{x}+\frac{1}{2}},y):=0,\label{eq:f_jump4}\\
\hat{H}_{z,h}^{n+\frac{1}{2}}(x,y_{j+\frac{1}{2}})
&:=H_{z,h}^{n+\frac{1}{2},-}(x,y_{j+\frac{1}{2}})\nonumber\\
&\quad\text{for all } j = 1, 2, 3,\ldots ,N_{y},\label{eq:f_jump5}\\
\hat{H}_{z,h}^{n+\frac{1}{2}}(x,y_{\frac{1}{2}})
&:=H_{z,h}^{n+\frac{1}{2},+}(x,y_{\frac{1}{2}})\nonumber\\
&+\frac{c_{0}}{2}\Big[ E_{x,h}^{n+1}(x,y_{\frac{1}{2}}^{+})
+E_{x,h}^{n}(x,y_{\frac{1}{2}}^{+}) \Big],\label{eq:f_jump6}\\
\hat{H}_{z,h}^{n+\frac{1}{2}}(x_{i+\frac{1}{2}},y)
&:=H_{z,h}^{n+\frac{1}{2},-}(x_{i+\frac{1}{2}},y)\nonumber\\
&\quad\text{for all } j = 1, 2, 3,\ldots ,N_{x},\label{eq:f_jump7}\\
\hat{H}_{z,h}^{n+\frac{1}{2}}(x_{\frac{1}{2}},y)
&:=H_{z,h}^{n+\frac{1}{2},+}(x_{\frac{1}{2}},y)\nonumber\\
&-\frac{c_{0}}{2}\Big[ E_{y,h}^{n+1}(x_{\frac{1}{2}}^{+},y)
+E_{x,h}^{n}(x_{\frac{1}{2}}^{+},y) \Big].\label{eq:f_jump8}
\end{align}
Due to the PEC condition \eqref{eq:pec_con} we have
$E_{x,h}^{n+1}(x,y_{\frac{1}{2}}^{+})
=E_{x,h}^{n+1}(x,y_{\frac{1}{2}}^{+})-E_{x,h}^{n+1}(x,y_{\frac{1}{2}}^{-})
=\llbracket E_{x,h}^{n+1}(x,y_{\frac{1}{2}}) \rrbracket$
in the equation \eqref{eq:f_jump6}, and the analogous one for
the other artificial viscosity in the equation \eqref{eq:f_jump8}.
The boundary terms are defined as follows:
\begin{align*}
\sigma_{Ih}:&=-\int_{I_{i}}[(\hat{H}_{z,h}^{n+\frac{1}{2}}(E_{x,h}^{n+1}
+E_{x,h}^{n})^{-})_{x,j+\frac{1}{2}}\\
&-(\hat{H}_{z,h}^{n+\frac{1}{2}}(E_{x,h}^{n+1}
+E_{x,h}^{n})^{+})_{x,j-\frac{1}{2}}]dx\\
&-\int_{I_{i}}[(\hat{E}_{x,h}^{n+1}(H_{z,h}^{n+\frac{3}{2}}
+H_{z,h}^{n+\frac{1}{2}})^{-})_{x,j+\frac{1}{2}}\\
&-(\hat{E}_{x,h}^{n+1}(H_{z,h}^{n+\frac{3}{2}}
+H_{z,h}^{n+\frac{1}{2}})^{+})_{x,j-\frac{1}{2}}]dx,
\\
\sigma_{Jh}:=&\int_{J_{j}}[(\hat{H}_{z,h}^{n+\frac{1}{2}}(E_{y,h}^{n+1}
+E_{y,h}^{n})^{-})_{i+\frac{1}{2},y}\\
&-(\hat{H}_{z,h}^{n+\frac{1}{2}}(E_{y,h}^{n+1}
+E_{y,h}^{n})^{+})_{i-\frac{1}{2},y}]dy\\
&+\int_{J_{j}}[(\hat{E}_{y,h}^{n+1}(H_{z,h}^{n+\frac{3}{2}}
+H_{z,h}^{n+\frac{1}{2}})^{-})_{i+\frac{1}{2},y}\\
&-(\hat{E}_{y,h}^{n+1}(H_{z,h}^{n+\frac{3}{2}}
+H_{z,h}^{n+\frac{1}{2}})^{+})_{i-\frac{1}{2},y}]dy.
\end{align*}
It should be noted that a nonlinear system of equations remains to be solved
in each time step. An investigation of nonlinear solvers, especially under
the aspect of energy conservation also for the approximations obtained with them,
is still pending.
However, we have had very positive experiences in the application of Newton
(or Newton-like) methods in solving such similar nonlinear problems
that arise when applying conforming methods \cite{Angermann:20a}. 

The proof of the energy relation in the subsequent section is based on
the following lemmas.
\begin{lemma}\label{lemma:fully_sum1}
For $n=1,2,\ldots,N$, with the flux densities \eqref{eq:f_jump1}--\eqref{eq:f_jump8}, we have
\begin{align*}
&\sum_{n=0}^{N}\sum_{i=1}^{N_{x}}\sum_{j=1}^{N_{y}}
\int_{K_{ij}}\Big[H_{z,h}^{n+\frac{1}{2}}\partial_{y}(E_{x,h}^{n+1}
+E_{x,h}^{n})\\
&+E_{x,h}^{n+1}\partial_{y}(H_{z,h}^{n+\frac{3}{2}}
+H_{z,h}^{n+\frac{1}{2}})\Big]
+\sum_{n=0}^{N}\sum_{i=1}^{N_{x}}\sum_{j=1}^{N_{y}}\sigma_{Ih}\\
&=\sum_{j=1}^{N_{y}-1}\Big[\int_{r}^{s}\Big(E_{x,h}^{N+1, +}
\llbracket H_{z,h}^{N+\frac{3}{2}}\rrbracket\Big)_{x,j+\frac{1}{2}}\\
&-\int_{r}^{s}\Big(E_{x,h}^{0, +}\llbracket H_{z,h}^{\frac{1}{2}}\rrbracket\Big)_{x,j+\frac{1}{2}} \Big]\\
&+\sum_{i=1}^{N_{x}}\sum_{j=1}^{N_{y}}\int_{K_{ij}}
\Big[E_{x,h}^{N+1}\partial_{y}H_{z,h}^{N+\frac{3}{2}}
-E_{x,h}^{0}\partial_{y}H_{z,h}^{\frac{1}{2}}\Big]\\
&+\frac{c_{0}}{2}\sum_{n=0}^{N}\int_{r}^{s}(E_{x,h}^{n+1,+}
+E_{x,h}^{n,+})^{2}_{x,\frac{1}{2}}.
\end{align*}
\end{lemma}
\proof
For details see \cite[Lemma 4.1]{Li:17b}.
\close
\begin{lemma}\label{lemma:fully_sum2}
For $n=1,2,\ldots,N$, with the flux densities \eqref{eq:f_jump1}--\eqref{eq:f_jump8}, we have
\begin{align*}
&-\sum_{n=0}^{N}\sum_{i=1}^{N_{x}}\sum_{j=1}^{N_{y}}\int_{K_{ij}}
\Big[H_{z,h}^{n+\frac{1}{2}}\partial_{x}(E_{y,h}^{n+1}+E_{y,h}^{n})\\
&+E_{y,h}^{n+1}\partial_{x}(H_{z,h}^{n+\frac{3}{2}}+H_{z,h}^{n+\frac{1}{2}})\Big]
+\sum_{n=0}^{N}\sum_{i=1}^{N_{x}}\sum_{j=1}^{N_{y}}\sigma_{Jh}\\
&=\sum_{i=1}^{N_{x}-1}\Big[-\int_{p}^{q}\Big(E_{y,h}^{N+1, +}
\llbracket H_{z,h}^{N+\frac{3}{2}}\rrbracket\Big)_{i+\frac{1}{2},y}\\
&+\int_{p}^{q}\Big(E_{y,h}^{0, +}\llbracket H_{z,h}^{\frac{1}{2}}\rrbracket\Big)_{i+\frac{1}{2},y}\Big]\\
&+\sum_{i=1}^{N_{x}}\sum_{j=1}^{N_{y}}\int_{K_{ij}}
\Big[-E_{y,h}^{N+1}\partial_{x}H_{z,h}^{N+\frac{3}{2}}
+E_{y,h}^{0}\partial_{x}H_{z,h}^{\frac{1}{2}}\Big]\\
&+\frac{c_{0}}{2}\sum_{n=0}^{N}\int_{p}^{q}(E_{y,h}^{n+1,+}
+E_{y,h}^{n,+})^{2}_{\frac{1}{2},y}.
\end{align*}
\end{lemma}
\proof
For details see \cite[Lemma 4.2]{Li:17b}.
\close

\bigskip
\section{The Nonlinear Electromagnetic Energy of the Full Discretization}
The nonlinear electromagnetic energy for the fully discrete approximation
(i.e.\ both in space and time) of the system \eqref{eq:f_x}--\eqref{eq:f_d2}
at $t^{n}$, $n=0,1,2,\ldots,N$, is defined by
\[
\mathcal{E}_{h}^{n}:=\|\EE_{h}^{n}\|_{\eps_{0}(1+\chi^{(1)})}^{2}
+\|H_{z,h}^{n+\frac{1}{2}}\|_{\mu_{0}}^{2}
+\big\|\vert\EE_h^{n}\vert^2\big\|_{\eps_{0}\chi^{(3)}}^{2}.
\]
In analogy to the conservativity and boundedness results for the continuous and
semi-discrete nonlinear electromagnetic energy
(Thms.~\ref{th:dg_energycontin}, \ref{th:dg_weaKenergy}),
in this section we demonstrate a stability result for
the fully discrete nonlinear electromagnetic energy
of the system \eqref{eq:f_x}--\eqref{eq:f_d2}.
\begin{theorem}\label{th:dg_fullenergy}
Let $(E_{x,h}^{n},E_{y,h}^{n},H_{z,h}^{n+\frac{1}{2}})^T\in (U_{h}^{k})^3$, $n\in\N$,
be the fully discrete solution of
\eqref{eq:f_x}--\eqref{eq:f_d2} for given $\JJ_h\in C(0,T,U_{h}^{k})^2$.
Then, if $\Delta t>0$, $h>0$ are sufficiently small and if $\Delta t/h$ is bounded
by some constant, the fully discrete nonlinear electromagnetic energy satisfies
\[
\mathcal{E}_{h}^{N}\le 3\mathcal{E}_{h}^{0}=3\Big[\|\EE_{h}^{0}\|_{\eps_{0}(1+\chi^{(1)})}^{2}
+\|H_{z,h}^{\frac{1}{2}}\|_{\mu_{0}}^{2}
+\big\|\vert\EE_h^{0}\vert^2\big\|_{\eps_{0}\chi^{(3)}}^{2}\Big]
\]
for vanishing current density
and
\[
\mathcal{E}_{h}^{N}
\le \exp(8T+1)\bigg[3\mathcal{E}_{h}^{0}
+\Delta t\,\sum_{n=0}^{N-1}\|\JJ_h^{n+\frac{1}{2}}\|_{(\eps_{0}(1+\chi^{(1)}))^{-1}}^{2}\bigg]
\]
for non-zero current density.
\end{theorem}
\proof
Taking $\Phi_{1h}:=(E_{x,h}^{n+1}+E_{x,h}^{n})$ in the equation \eqref{eq:f_d1},
we have
\begin{align}
&\int_{K_{ij}}(D_{x,h}^{n+1}-D_{x,h}^{n})(E_{x,h}^{n+1}+E_{x,h}^{n})\nonumber\\
&=\int_{K_{ij}}\eps_{0}(1+\chi^{(1)})(E_{x,h}^{n+1}-E_{x,h}^{n})(E_{x,h}^{n+1}+E_{x,h}^{n
})\nonumber\\
&\quad + \int_{K_{ij}}\eps_{0}\chi^{(3)}\Big[\frac{1}{2}\big[(E_{x,h}^{
n+1})^{2}+(E_{x,h}^{n})^{2}+(E_{y,h}^{n+1})^{2}+(E_{y,h}^{n})^{2}\big]\nonumber\\
&\quad \times
(E_{x,h}^{n+1}-E_{x,h}^{n})(E_{x,h}^{n+1}+E_{x,h}^{n})\nonumber\\
&\quad +\Big([(E_{x,h}^{n+1})^{2}+(E_{x,h}^{n})^{2}](E_{x,h}^{n+1}
-E_{x,h}^{n})(E_{x,h}^{n+1}+E_{x,h}^{n}
)\nonumber\\
&\quad +[E_{x,h}^{n+1}E_{y,h}^{n+1}+E_{x,h}^{n}E_{y,h}^{n}](E_{y,h}^{n+1}
-E_{y,h}^{n})(E_{x,h}^{n+1}+E_{x,h}^{n})\Big)\Big].\label{eq:f_et1}
\end{align}
Taking $\Phi_{2h}:=(E_{y,h}^{n+1}+E_{y,h}^{n})$ in the equation \eqref{eq:f_d2},
we have
\begin{align}
&\int_{K_{ij}}(D_{y,h}^{n+1}-D_{y,h}^{n})(E_{y,h}^{n+1}+E_{y,h}^{n})\nonumber\\
&=\int_{K_{ij}}\eps_{0}(1+\chi^{(1)})(E_{y,h}^{n+1}-E_{y,h}^{n})(E_{y,h}^{n+1}+E_{y,h}^{n
})\nonumber\\
&\quad +\int_{K_{ij}}\eps_{0}\chi^{(3)}\Big[\frac{1}{2}\big[(E_{x,h}^{
n+1})^{2}+(E_{x,h}^{n})^{2}+(E_{y,h}^{n+1})^{2}+(E_{y,h}^{n})^{2}\big]\nonumber\\
&\quad \times
(E_{y,h}^{n+1}-E_{y,h}^{n})(E_{y,h}^{n+1}+E_{y,h}^{n})\nonumber\\
&\quad +\Big([(E_{y,h}^{n+1})^{2}+(E_{y,h}^{n})^{2}](E_{y,h}^{n+1}
-E_{y,h}^{n})(E_{y,h}^{n+1}+E_{y,h}^{n}
)\nonumber\\
&\quad +[E_{x,h}^{n+1}E_{y,h}^{n+1}+E_{x,h}^{n}E_{y,h}^{n}](E_{x,h}^{n+1}
-E_{x,h}^{n})(E_{y,h}^{n+1}+E_{y,h}^{n})\Big)\Big].\label{eq:f_et2}
\end{align}
Adding the equations \eqref{eq:f_et1} and \eqref{eq:f_et2}, we see that
\begin{align*}
&\int_{K_{ij}}(D_{x,h}^{n+1}-D_{x,h}^{n})(E_{x,h}^{n+1}+E_{x,h}^{n})\\
&+\int_{K_{i,j}}(D_{y,h}^{n+1}-D_{y,h}^{n})(E_{y,h}^{n+1}+E_{y,h}^{n})\\
&=\int_{K_{ij}}\eps_{0}(1+\chi^{(1)})
\big[\vert\EE_h^{n+1}\vert^2-\vert\EE_h^n\vert^2\big]
\\
&\quad +\int_{K_{ij}}\eps_{0}\chi^{(3)}\Big[\frac{1}{2}
\big[\vert\EE_h^{n+1}\vert^2+\vert\EE_h^n\vert^2\big]\big[\vert\EE_h^{n+1}\vert^2-\vert\EE_h^n\vert^2\big]
\\
&\quad +[(E_{x,h}^{n+1})^{2}+(E_{x,h}^{n})^{2}][(E_{x,h}^{n+1})^{2}-(E_{x,h}^{n})^{2}]\\
&\quad +[(E_{y,h}^{n+1})^{2}+(E_{y,h}^{n})^{2}]
[(E_{y,h}^{n+1})^{2}-(E_{y,h}^{n})^{2}]\\
&\quad +[E_{x,h}^{n+1}E_{y,h}^{n+1}+E_{x,h}^{n}E_{y,h}^{n}][E_{x,h}^{n+1}E_{y,h}^{n+1}
+E_{x,h}^{n}E_{y,h}^{n+1}\\
&\quad -E_{x,h}^{n+1}E_{y,h}^{n}-E_{y,h}^{n}E_{x,h}^{n}]\\
&\quad +[E_{x,h}^{n+1}E_{y,h}^{n+1}+E_{x,h}^{n}E_{y,h}^{n}][E_{x,h}^{n+1}E_{y,h}^{n+1}
+E_{x,h}^{n+1}E_{y,h}^{n}\\
&\quad -E_{x,h}^{n}E_{y,h}^{n+1}-E_{x,h}^{n}E_{y,h}^{n}]\Big].
\end{align*}
The term in square brackets in the second integral of the right-hand side
can be simplified as follows:
\begin{align*}
\Big[\ldots\Big]
&=\frac{1}{2}
\big[\vert\EE_h^{n+1}\vert^2+\vert\EE_h^n\vert^2\big]\big[\vert\EE_h^{n+1}\vert^2-\vert\EE_h^n\vert^2\big]\\
&\quad + \vert E_{x,h}^{n+1}\vert^4-\vert E_{x,h}^{n}\vert^4 + \vert E_{y,h}^{n+1}\vert^4-\vert E_{y,h}^{n}\vert^4\\
&\quad + 2\big[E_{x,h}^{n+1}E_{y,h}^{n+1}+E_{x,h}^{n}E_{y,h}^{n}\big]
\big[E_{x,h}^{n+1}E_{y,h}^{n+1}-E_{x,h}^{n}E_{y,h}^{n}\big]\\
&=\frac{1}{2}\big[\vert\EE_h^{n+1}\vert^4-\vert\EE_h^n\vert^4\big]\\
&\quad + \vert E_{x,h}^{n+1}\vert^4 + \vert E_{y,h}^{n+1}\vert^4
- \vert E_{x,h}^{n}\vert^4 - \vert E_{y,h}^{n}\vert^4\\
&+ 2\vert E_{x,h}^{n+1}\vert^2\vert E_{y,h}^{n+1}\vert^2 - 2\vert E_{x,h}^{n}\vert^2\vert E_{y,h}^{n}\vert^2\\
&=\frac{3}{2}\big[\vert\EE_h^{n+1}\vert^4-\vert\EE_h^n\vert^4\big].
\end{align*}
Hence we get
\begin{align}
&\int_{K_{ij}}(D_{x,h}^{n+1}-D_{x,h}^{n})(E_{x,h}^{n+1}+E_{x,h}^{n})\nonumber\\
&+\int_{K_{ij}}(D_{y,h}^{n+1}-D_{y,h}^{n})(E_{y,h}^{n+1}+E_{y,h}^{n})\nonumber\\
&=\int_{K_{ij}}\eps_{0}(1+\chi^{(1)})\big[\vert\EE_h^{n+1}\vert^2-\vert\EE_h^n\vert^2\big]\nonumber\\
&+\frac{3}{2}\int_{K_{ij}}\eps_{0}\chi^{(3)}\big[\vert\EE_h^{n+1}\vert^4-\vert\EE_h^n\vert^4\big].
\label{eq:f_et6}
\end{align}
Taking $\Phi_{1h}:=2\Delta t(E_{x,h}^{n+1}+E_{x,h}^{n})$ in the equation
\eqref{eq:f_x}, we have
\begin{align}
&2\int_{K_{ij}}(D_{x,h}^{n+1}-D_{x,h}^{n})(E_{x,h}^{n+1}+E_{x,h}^{n})\nonumber\\
&-2\Delta t\int_{I_{i}}[(\hat{H}_{z,h}^{n+\frac{1}{2}}(E_{x,h}^{n+1}
+E_{x,h}^{n})^{-})_{x,j+\frac{1}{2}}\nonumber\\
&-(\hat{H}_{z,h}^{n+\frac{1}{2}}(E_{x,h}^{n+1}+E_{x,h}^{n})^{+})_{x,j-\frac{1}{2}}]dx\nonumber\\
&+2\Delta t\int_{K_{ij}}H_{z,h}^{n+\frac{1}{2}}\partial_{y}(E_{x,h}^{n+1}+E_{x,h}^{n})\nonumber\\
&-2\Delta t\int_{K_{ij}}J_{x,h}^{n+\frac{1}{2}}(E_{x,h}^{n+1}+E_{x,h}^{n})=0.
\label{eq:f_et7}
\end{align}
Taking $\Phi_{2h}:=2\Delta t(E_{y,h}^{n+1}+E_{y,h}^{n})$ in the equation
\eqref{eq:f_y}, we have
\begin{align}
&2\int_{K_{ij}}(D_{y,h}^{n+1}-D_{y,h)}(E_{y,h}^{n+1}+E_{y,h}^{n})\nonumber\\
&+2\Delta t\int_{J_{j}}[(\hat{H}_{z,h}^{n+\frac{1}{2}}(E_{y,h}^{n+1}
+E_{y,h}^{n})^{-})_{i+\frac{1}{2},y}\nonumber\\
&-(\hat{H}_{z,h}^{n+\frac{1}{2}}(E_{y,h}^{n+1}+E_{y,h
}^{n})^{+})_{i-\frac{1}{2},y}]dy\nonumber\\
&-2\Delta t\int_{K_{ij}}H_{z,h}^{n+\frac{1}{2}}\partial_{x}(E_{y,h}^{n+1}+E_{y,h}^{n})\nonumber\\
&-2\Delta
t\int_{K_{ij}}J_{y,h}^{n+\frac{1}{2}}(E_{y,h}^{n+1}+E_{y,h}^{n})=0.
\label{eq:f_et8}
\end{align}
Taking $\Phi_{3h}:=2\Delta t(H_{z,h}^{n+\frac{3}{2}}+H_{z,h}^{n+\frac{1}{2}})$
in the equation \eqref{eq:f_h}, we have
\begin{align*}
&2\int_{K_{ij}}\mu_{0}(H_{z,h}^{n+\frac{3}{2}}-H_{z,h}^{n+\frac{1}{2}})(H_{z,h}
^{n+\frac{3}{2}}+H_{z,h}^{n+\frac{1}{2}})\\
&+2\Delta t\int_{J_{j}}[(\hat{E}_{y,h}^{n+1}(H_{z,h}^{n+\frac{3}{2}}+H_{z,h}^{n+\frac{1}{2
}})^{-})_{i+\frac{1}{2},y}\\
&-(\hat{E}_{y,h}^{n+1}(H_{z,h}^{n+\frac{3}{2}}+H_{z,h}^{n+\frac{1}{2}})^{+})_{i-\frac{1}{2},y}]dy\\
&-2\Delta t\int_{K_{ij}}E_{y,h}^{n+1}\partial_{x}(H_{z,h}^{n+\frac{3}{2}}+H_{z,h}^{n+\frac{1}{2}})\\
&-2\Delta t\int_{I_{i}}[(\hat{E}_{x,h}^{n+1}(H_{z,h}^{n+\frac{3}{2}}+H_{z,h}^{n+\frac{1}{2
}})^{-})_{x,j+\frac{1}{2}}\\
&-(\hat{E}_{x,h}^{n+1}(H_{z,h}^{n+\frac{3}{2}}+H_{z,h}^
{n+\frac{1}{2}})^{+})_{x,j-\frac{1}{2}}]dx\\
&+2\Delta t\int_{K_{ij}}E_{x,h}^{n+1}\partial_{y}(H_{z,h}^{n+\frac{3}{2}}+H_{z,h}^{
n+\frac{1}{2}})=0.
\end{align*}
Adding the equations \eqref{eq:f_et7} and \eqref{eq:f_et8},
substituting the result in equation \eqref{eq:f_et6}, we obtain
\begin{align*}
&2\int_{K_{ij}}\eps_{0}(1+\chi^{(1)})\big[\vert\EE_h^{n+1}\vert^2-\vert\EE_h^n\vert^2\big]\\
&+2\int_{K_{ij}}\mu_{0}[(H_{z,h}^{n+\frac{3}{2}})^{2}-(H_{z,h}^{n+\frac{1}{2}})^{2}]\\
&+3\int_{K_{ij}}\eps_{0}\chi^{(3)}\big[\vert\EE_h^{n+1}\vert^4-\vert\EE_h^n\vert^4\big]\\
&-2\Delta t\int_{I_{i}}[(\hat{H}_{z,h}^{n+\frac{1}{2}}(E_{x,h}^{n+1}
+E_{x,h}^{n})^{-})_{x,j+\frac{1}{2}}\\
&-(\hat{H}_{z,h}^{n+\frac{1}{2}}(E_{x,h}^{n+1}+E_{x,h}^{n})^{+})_{x,j-\frac{1}{2}}]dx\\
&+2\Delta t\int_{K_{ij}}H_{z,h}^{n+\frac{1}{2}}\partial_{y}(E_{x,h}^{n+1}+E_{x,h}^{n})\\
&+2\Delta t\int_{J_{j}}[(\hat{H}_{z,h}^{n+\frac{1}{2}}(E_{y,h}^{n+1}
+E_{y,h}^{n})^{-})_{i+\frac{1}{2},y}\\
&-(\hat{H}_{z,h}^{n+\frac{1}{2}}(E_{y,h}^{n+1}
+E_{y,h}^{n})^{+})_{i-\frac{1}{2},y}]dy\\
&-2\Delta t\int_{K_{ij}}H_{z,h}^{n+\frac{1}{2}}\partial_{x}(E_{y,h}^{n+1}+E_{y,h}^{n})\\
&+2\Delta t\int_{J_{j}}[(\hat{E}_{y,h}^{n+1}(H_{z,h}^{n+\frac{3}{2}}+H_{z,h}^{n+\frac{1}{2}})^{-})_{i+\frac{1}{2},y}\\
&-(\hat{E}_{y,h}^{n+1}(H_{z,h}^{n+\frac{3}{2}}
+H_{z,h}^{n+\frac{1}{2}})^{+})_{i-\frac{1}{2},y}]dy\\
&-2\Delta t\int_{K_{ij}}E_{y,h}^{n+1}\partial_{x}(H_{z,h}^{n+\frac{3}{2}}+H_{z,h}^{n+\frac{1}{2}})\\
&-2\Delta t\int_{I_{i}}[(\hat{E}_{x,h}^{n+1}(H_{z,h}^{n+\frac{3}{2}}
+H_{z,h}^{n+\frac{1}{2}})^{-})_{x,j+\frac{1}{2}}\\
&-(\hat{E}_{x,h}^{n+1}(H_{z,h}^{n+\frac{3}{2}}+H_{z,h}^{n+\frac{1}{2}})^{+})_{x,j-\frac{1}{2}}]dx\\
&+2\Delta t\int_{K_{ij}}E_{x,h}^{n+1}\partial_{y}(H_{z,h}^{n+\frac{3}{2}}+H_{z,h}^{n+\frac{1}{2}})\\
&=2\Delta t\int_{K_{ij}}J_{x,h}^{n+\frac{1}{2}}(E_{x,h}^{n+1}+E_{x,h}^{n})\\
&+2\Delta t\int_{K_{ij}}J_{y,h}^{n+\frac{1}{2}}(E_{y,h}^{n+1}+E_{y,h}^{n}).
\end{align*}
Summing up over the $1\le i\le N_{x}$, $1\le i\le N_{y}$, and with respect
to time from $n=1$ to $N$, and using the Lemmas
\ref{lemma:fully_sum1}--\ref{lemma:fully_sum2} we arrive at
\begin{align}
&2\Big[\|\EE_h^{N+1}\|_{\eps_{0}(1+\chi^{(1)})}^{2}
-\|\EE_h^{0}\|_{\eps_{0}(1+\chi^{(1)})}^{2}\nonumber\\
&+\|H_{z,h}^{N+\frac{3}{2}}\|_{\mu_{0}}^{2}-\|H_{z,h}^{\frac{1}{2}}\|_{\mu_{0}}^{2}\Big]\nonumber\\
&+3\Big[\|\vert\EE_h^{N+1}\vert^2\|_{\eps_{0}\chi^{(3)}}^2
-\|\vert\EE_h^{0}\vert^2\|_{\eps_{0}\chi^{(3)}}^2\Big]\nonumber\\
&=2\Delta t\sum_{n=0}^{N}\int_{\Omega}J_{y,h}^{n+\frac{1}{2}}(E_{y,h}^{n+1}+E_{y,h}^{n})\nonumber\\
&+2\Delta t\sum_{n=0}^{N}\int_{\Omega}J_{x,h}^{n+\frac{1}{2}}(E_{x,h}^{n+1}+E_{x,h}^{n})\nonumber\\
&-2B_{y}\Big(E_{x,h}^{N+1}, H_{z,h}^{N+\frac{3}{2}}\Big)
+2B_{y}\Big(E_{x,h}^{0}, H_{z,h}^{\frac{1}{2}}\Big)\nonumber\\
&+2B_{x}\Big(E_{y,h}^{N+1}, H_{z,h}^{N+\frac{3}{2}}\Big)
-2B_{x}\Big(E_{y,h}^{0}, H_{z,h}^{\frac{1}{2}}\Big),
\label{eq:f_et11}
\end{align}
where the bilinear forms are defined as
\begin{align*}
B_{x}\Big(E_{y,h}^{n+1}, H_{z,h}^{n+\frac{3}{2}}\Big)
&:=\Delta t \Big[\sum_{i=1}^{N_{x}}\sum_{j=1}^{N_{y}}
\int_{K_{ij}}E_{y,h}^{n+1}\partial_{x}H_{z,h}^{n+\frac{3}{2}}
\\&\quad
+ \sum_{i=1}^{N_{x}-1}\int_{p}^{q}(E_{y,h})_{i+\frac{1}{2}}^{+}
\llbracket H_{z,h}^{n+\frac{3}{2}}\rrbracket_{i+\frac{1}{2}}dy\Big],
\\
B_{y}\Big(E_{x,h}^{n+1}, H_{z,h}^{n+\frac{3}{2}}\Big)
&:=\Delta t \Big[\sum_{i=1}^{N_{x}}\sum_{j=1}^{N_{y}}
\int_{K_{ij}}E_{x,h}^{n+1}\partial_{y}H_{z,h}^{n+\frac{3}{2}}
\\&\quad
+ \sum_{j=1}^{N_{y}-1}\int_{r}^{s}(E_{x,h})_{j+\frac{1}{2}}^{+}
\llbracket H_{z,h}^{n+\frac{3}{2}}\rrbracket_{j+\frac{1}{2}}dx \Big]
\end{align*}
(cf.\ \cite[Proof of Thm.~4.1]{Li:17b} or \cite[eq.~(4.1)]{Zhang:10b}).
Using an inverse estimate (cf.\ \cite[Proof of Thm.~4.1]{Li:17b} or \cite[Lemma 4.1]{Zhang:10b}),
we have that
\begin{align*}
B_{y}\Big(E_{x,h}^{n+1}, H_{z,h}^{n+\frac{3}{2}}\Big)
\le 2 \Delta t C_{INV}\frac{C_{\eps\mu}}{h}
\|E_{x,h}^{n+1}\|_{\eps_{0}(1+\chi^{(1)})}\|H_{z,h}^{n+\frac{3}{2}}\|_{\mu_{0}},
\end{align*}
where $C_{INV}$ is a positive constant that is independent of $h$ and $\Delta t$,
and $C_{\eps\mu}:=\big\|(\eps_{0}\mu_{0}(1+\chi^{(1)}))^{-1/2}\big\|_{L_\infty(\Omega)}$.
The right-hand side is estimated by means of Young's inequality with $\eps$
(see, e.g., \cite[Lemma 1, 2)]{Angermann:19g}),
where the parameter called here $\alpha>0$ will be determined later:
\begin{equation}\label{eq:f_et14}
\begin{aligned}
B_{y}\Big(E_{x,h}^{n+1}, H_{z,h}^{n+\frac{3}{2}}\Big)
\le \alpha\|E_{x,h}^{n+1}\|_{\eps_{0}(1+\chi^{(1)})}^{2}
+\Big(\Delta t C_{INV}\frac{C_{\eps\mu}}{\alpha h}\Big)^2
\|H_{z,h}^{n+\frac{3}{2}}\|_{\mu_{0}}^{2}.
\end{aligned}
\end{equation}
Similarly we get (with the same parameter $\alpha$)
\begin{equation}\label{eq:f_et15}
\begin{aligned}
B_{x}\Big(E_{x,y}^{n+1}, H_{z,h}^{n+\frac{3}{2}}\Big)
\le \alpha\|E_{y,h}^{n+1}\|_{\eps_{0}(1+\chi^{(1)})}^{2}
+\Big(\Delta t C_{INV}\frac{C_{\eps\mu}}{\alpha h}\Big)^2
\|H_{z,h}^{n+\frac{3}{2}}\|_{\mu_{0}}^{2},
\end{aligned}
\end{equation}
and
\begin{align}
B_{y}\Big(E_{x,h}^{0}, H_{z,h}^{\frac{1}{2}}\Big)
&
\le \Delta t C_{INV}\frac{C_{\eps\mu}}{h}\Big[\|E_{x,h}^{0}\|_{\eps_{0}(1+\chi^{(1)})}^{2}
+\|H_{z,h}^{\frac{1}{2}}\|_{\mu_{0}}^{2}\Big],\label{eq:f_et16}\\
B_{x}\Big(E_{y,h}^{0}, H_{z,h}^{\frac{1}{2}}\Big)
&\le \Delta t C_{INV}\frac{C_{\eps\mu}}{h}\Big[\|E_{y,h}^{0}\|_{\eps_{0}(1+\chi^{(1)})}^{2}
+\|H_{z,h}^{\frac{1}{2}}\|_{\mu_{0}}^{2}\Big].\label{eq:f_et17}
\end{align}
The first two terms from the right-hand side of equation \eqref{eq:f_et11}
are estimated by means of Young's inequality, too. This gives
\begin{align}
&2\Delta t\,\sum_{n=0}^{N}\int_{\Omega}J_{x,h}^{n+\frac{1}{2}}(E_{x,h}^{n+1}+E_{x,h}^{n})\nonumber\\
&=2\Delta t\,\sum_{n=0}^{N}\int_{\Omega}
(\eps_{0}(1+\chi^{(1)}))^{-1/2}J_{x,h}^{n+\frac{1}{2}}
((\eps_{0}(1+\chi^{(1)}))^{1/2}(E_{x,h}^{n+1}+E_{x,h}^{n})\nonumber\\
&\le \Delta t\,\sum_{n=0}^{N}\|J_{x,h}^{n+\frac{1}{2}}\|_{(\eps_{0}(1+\chi^{(1)}))^{-1}}^{2}
+ \Delta t\,\sum_{n=0}^{N}\|E_{x,h}^{n+1}+E_{x,h}^{n}\|_{\eps_{0}(1+\chi^{(1)})}^{2}
\nonumber\\
&\le \Delta t\,\sum_{n=0}^{N}\|J_{x,h}^{n+\frac{1}{2}}\|_{(\eps_{0}(1+\chi^{(1)}))^{-1}}^{2}\nonumber\\
&+ 2\Delta t\,\sum_{n=0}^{N}\Big[\|E_{x,h}^{n+1}\|_{\eps_{0}(1+\chi^{(1)})}^{2}
+\|E_{x,h}^{n}\|_{\eps_{0}(1+\chi^{(1)})}^{2}\Big]\nonumber\\
&\le \Delta t\,\sum_{n=0}^{N}\|J_{x,h}^{n+\frac{1}{2}}\|_{(\eps_{0}(1+\chi^{(1)}))^{-1}}^{2}\nonumber\\
&+ 4\Delta t\,\sum_{n=0}^{N}\|E_{x,h}^{n+1}\|_{\eps_{0}(1+\chi^{(1)})}^{2}
+ 2\Delta t\,\|E_{x,h}^{0}\|_{\eps_{0}(1+\chi^{(1)})}^{2}
\label{eq:f_et12}
\end{align}
and
\begin{align}
&2\Delta t\,\sum_{n=0}^{N}\int_{\Omega}J_{y,h}^{n+\frac{1}{2}}(E_{y,h}^{n+1}+E_{y,h}^{n})\nonumber\\
&\le \Delta t\,\sum_{n=0}^{N}\|J_{y,h}^{n+\frac{1}{2}}\|_{(\eps_{0}(1+\chi^{(1)}))^{-1}}^{2}\nonumber\\
&+ 4\Delta t\,\sum_{n=0}^{N}\|E_{y,h}^{n+1}\|_{\eps_{0}(1+\chi^{(1)})}^{2}
+ 2\Delta t\,\|E_{y,h}^{0}\|_{\eps_{0}(1+\chi^{(1)})}^{2}.\label{eq:f_et13}
\end{align}
Finally, using the estimates \eqref{eq:f_et12}, \eqref{eq:f_et13}, and
\eqref{eq:f_et14}--\eqref{eq:f_et17} in \eqref{eq:f_et11}, we obtain
\begin{align*}
&2\Big[\|\EE_h^{N+1}\|_{\eps_{0}(1+\chi^{(1)})}^{2}+\|H_{z,h}^{N+\frac{3}{2}}\|_{\mu_{0}}^{2}\Big]
+\|\vert\EE_h^{N+1}\vert^2\|_{\eps_{0}\chi^{(3)}}^2\\
&\le 4\Delta t\,\sum_{n=0}^{N}\|\EE_h^{n+1}\|_{\eps_{0}(1+\chi^{(1)})}^{2}
+ 2\Delta t\,\|\EE_h^{0}\|_{\eps_{0}(1+\chi^{(1)})}^{2}\\
& + 2\alpha\|\EE_h^{N+1}\|_{\eps_{0}(1+\chi^{(1)})}^{2}
+4\Big(\Delta t\, C_{INV}\frac{C_{\eps\mu}}{\alpha h}\Big)^2
\|H_{z,h}^{N+\frac{3}{2}}\|_{\mu_{0}}^{2}\\
&+\Delta t\,\sum_{n=0}^{N}\|\JJ_h^{n+\frac{1}{2}}\|_{(\eps_{0}(1+\chi^{(1)}))^{-1}}^{2}\\
&+ 2\Delta t\, C_{INV}\frac{C_{\eps\mu}}{h}\Big[\|\EE_h^{0}\|_{\eps_{0}(1+\chi^{(1)})}^{2}
+2\|H_{z,h}^{\frac{1}{2}}\|_{\mu_{0}}^{2}\Big]\\
&+2\Big[\|\EE_h^{0}\|_{\eps_{0}(1+\chi^{(1)})}^{2}+\|H_{z,h}^{\frac{1}{2}}\|_{\mu_{0}}^{2}\Big]
+3\|\vert\EE_h^{0}\vert^2\|_{\eps_{0}\chi^{(3)}}^2.
\end{align*}
Now we chose $\alpha:=1/2$ and move the corresponding term to the left-hand side.
If the condition
\begin{equation}\label{eq:stepsizecond1}
\frac{\Delta t}{h}\le \min\left\{\frac{1}{4C_{INV}C_{\eps\mu}};\frac{1}{4h}\right\}
\end{equation}
is satisfied, we obtain
\begin{align*}
&\|\EE_h^{N+1}\|_{\eps_{0}(1+\chi^{(1)})}^{2}+\|H_{z,h}^{N+\frac{3}{2}}\|_{\mu_{0}}^{2}
+\|\vert\EE_h^{N+1}\vert^2\|_{\eps_{0}\chi^{(3)}}^2\\
&\le 4\Delta t\,\sum_{n=0}^{N}\|\EE_h^{n+1}\|_{\eps_{0}(1+\chi^{(1)})}^{2}
+ 2\Delta t\,\|\EE_h^{0}\|_{\eps_{0}(1+\chi^{(1)})}^{2}
+ \Delta t\,\sum_{n=0}^{N}\|\JJ_h^{n+\frac{1}{2}}\|_{(\eps_{0}(1+\chi^{(1)}))^{-1}}^{2}\\
&+ \frac{1}{2}\Big[\|\EE_h^{0}\|_{\eps_{0}(1+\chi^{(1)})}^{2}
+ 2\|H_{z,h}^{\frac{1}{2}}\|_{\mu_{0}}^{2}\Big]\\
&+2\Big[\|\EE_h^{0}\|_{\eps_{0}(1+\chi^{(1)})}^{2}+\|H_{z,h}^{\frac{1}{2}}\|_{\mu_{0}}^{2}\Big]
+ 3\|\vert\EE_h^{0}\vert^2\|_{\eps_{0}\chi^{(3)}}^2\\
&\le 4\Delta t\,\sum_{n=0}^{N}\|\EE_h^{n+1}\|_{\eps_{0}(1+\chi^{(1)})}^{2}
+ \Delta t\,\sum_{n=0}^{N}\|\JJ_h^{n+\frac{1}{2}}\|_{(\eps_{0}(1+\chi^{(1)}))^{-1}}^{2}\\
&+ 3\Big[\|\EE_h^{0}\|_{\eps_{0}(1+\chi^{(1)})}^{2}
+ \|H_{z,h}^{\frac{1}{2}}\|_{\mu_{0}}^{2}
+ \|\vert\EE_h^{0}\vert^2\|_{\eps_{0}\chi^{(3)}}^2\Big].
\end{align*}
If we strengthen the condition \eqref{eq:stepsizecond1} to
\[
\frac{\Delta t}{h} < \min\left\{\frac{1}{4C_{INV}C_{\eps\mu}};\frac{1}{4h}\right\},
\]
then we may apply a discrete Gronwall's inequality \cite[Lemma 5.1]{Heywood:90}
(also cited in \cite[Lemma 2]{Angermann:19g})
to obtain
\begin{align*}
&\|\EE_h^{N+1}\|_{\eps_{0}(1+\chi^{(1)})}^{2}+\|H_{z,h}^{N+\frac{3}{2}}\|_{\mu_{0}}^{2}
+\|\vert\EE_h^{N+1}\vert^2\|_{\eps_{0}\chi^{(3)}}^2\\
&\le \exp\bigg(4\Delta t \sum_{n=0}^{N}(1-4\Delta t)^{-1}\bigg)
\bigg[\Delta t\,\sum_{n=0}^{N}\|\JJ_h^{n+\frac{1}{2}}\|_{(\eps_{0}(1+\chi^{(1)}))^{-1}}^{2}\\
&+ 3\Big[\|\EE_h^{0}\|_{\eps_{0}(1+\chi^{(1)})}^{2}
+ \|H_{z,h}^{\frac{1}{2}}\|_{\mu_{0}}^{2}
+ \|\vert\EE_h^{0}\vert^2\|_{\eps_{0}\chi^{(3)}}^2\Big]\bigg].
\end{align*}
If even
\[
\frac{\Delta t}{h} \le \min\left\{\frac{1}{4C_{INV}C_{\eps\mu}};\frac{1}{8h}\right\},
\]
holds, then
\begin{align*}
&\|\EE_h^{N+1}\|_{\eps_{0}(1+\chi^{(1)})}^{2}+\|H_{z,h}^{N+\frac{3}{2}}\|_{\mu_{0}}^{2}
+\|\vert\EE_h^{N+1}\vert^2\|_{\eps_{0}\chi^{(3)}}^2\\
&\le \exp(8T+1)
\bigg[\Delta t\,\sum_{n=0}^{N}\|\JJ_h^{n+\frac{1}{2}}\|_{(\eps_{0}(1+\chi^{(1)}))^{-1}}^{2}\\
&+ 3\Big[\|\EE_h^{0}\|_{\eps_{0}(1+\chi^{(1)})}^{2}
+ \|H_{z,h}^{\frac{1}{2}}\|_{\mu_{0}}^{2}
+ \|\vert\EE_h^{0}\vert^2\|_{\eps_{0}\chi^{(3)}}^2\Big]\bigg].
\end{align*}
Since the term
$\Delta t\,\sum_{n=0}^{N}\|\JJ_h^{n+\frac{1}{2}}\|_{(\eps_{0}(1+\chi^{(1)}))^{-1}}^{2}$
can be interpreted as an approximation to
\linebreak
$\int_{0}^{T}\|\JJ(s)\|_{(\eps_{0}(1+\chi^{(1)}))^{-1}}^{2}ds$
it can be regarded as being bounded independently of $h$.

To prove the first statement, the estimates \eqref{eq:f_et12}, \eqref{eq:f_et13}
are not needed, and we immediately get from \eqref{eq:f_et11} the relation
\begin{align*}
&2\Big[\|\EE_h^{N+1}\|_{\eps_{0}(1+\chi^{(1)})}^{2}+\|H_{z,h}^{N+\frac{3}{2}}\|_{\mu_{0}}^{2}\Big]
+\|\vert\EE_h^{N+1}\vert^2\|_{\eps_{0}\chi^{(3)}}^2\\
&\le 2\alpha\|\EE_h^{N+1}\|_{\eps_{0}(1+\chi^{(1)})}^{2}
+4\Big(\Delta t\, C_{INV}\frac{C_{\eps\mu}}{\alpha h}\Big)^2
\|H_{z,h}^{N+\frac{3}{2}}\|_{\mu_{0}}^{2}\\
&+ 2\Delta t\, C_{INV}\frac{C_{\eps\mu}}{h}\Big[\|\EE_h^{0}\|_{\eps_{0}(1+\chi^{(1)})}^{2}
+2\|H_{z,h}^{\frac{1}{2}}\|_{\mu_{0}}^{2}\Big]\\
&+2\Big[\|\EE_h^{0}\|_{\eps_{0}(1+\chi^{(1)})}^{2}+\|H_{z,h}^{\frac{1}{2}}\|_{\mu_{0}}^{2}\Big]
+3\|\vert\EE_h^{0}\vert^2\|_{\eps_{0}\chi^{(3)}}^2.
\end{align*}
Now the condition
\[
\frac{\Delta t}{h}\le \frac{1}{4C_{INV}C_{\eps\mu}}
\]
already leads to the statement.
\close

\bigskip
\section{Error behavior of the Fully Discrete Solution}

If the assumptions of Thm.~\ref{th:dg_semi_error_estimate}
and Thm.~\ref{th:dg_fullenergy} are combined with the additional requirements
that the weak solution $(E_{x},E_{y},H_{z})^T$
of the system \eqref{eq:2dsystem} belongs to $C^2(0,T,H^{k+1}(\Omega))^3$, $k\in\N$,
the fully discrete solution
$(E_{x,h}^{n},E_{y,h}^{n},H_{z,h}^{n+\frac{1}{2}})^T\in (U_{h}^{k})^3$
of \eqref{eq:f_x}--\eqref{eq:f_d2} is uniformly bounded w.r.t.\ $h$ and $n\in\N$
and the initial values are chosen such that
$\mathcal{E}_{h}^{0}\le Ch^{2(k+1)}$ is satisfied,
then it is possible to prove a bound for the norm
\[
\|\EE_h^{N}-\EE(T)\|_{\eps_{0}(1+\chi^{(1)})}
+\|H_{z,h}^{N+\frac{1}{2}}-H_z(T)\|_{\mu_{0}}
\]
of the error of optimal order,
i.e.\ of the type $C(h^{k+1} + (\Delta t)^2)$.

The proof is based on the stability result Thm.~\ref{th:dg_fullenergy} and
runs structurally like the proof of Thm.~\ref{th:dg_semi_error_estimate},
whereby on the one hand the assumed boundedness of the fully discrete solution
(similar to the proof of Thm.~\ref{th:dg_semi_error_estimate})
and on the other hand standard estimates for time discretizations
(cf.\ \cite[Sect.~9.8]{Angermann:21a}) are used.

We do not want to describe the proof in detail, not only because it is quite technical
(and therefore very lengthy), but above all because we see a conceptual discrepancy
between the fact that on the one hand the introduced family of spatial dG discretizations
can be shown to be energy stable (see Thm.~\ref{th:dg_energycontin}),
while on the other hand -- as far as known to the authors -- (nonlinear) results analogous to
Thm.~\ref{th:dg_fullenergy} are only available for a few selected temporal discretization methods
of first and second order.
Although there is active research on methods that are aimed at establishing or improving
certain conservation properties (for instance implicit Runge-Kutta methods \cite{Hochbruck:15},
implicit-explicit Runge-Kutta (IMEX-RK) methods \cite{Boscarino:17} with
an appropiately chosen IMEX strategy, or symplectic methods \cite{sha:08}),
most of the theoretical results (if any) are related to the classical (linear) Maxwell system.
To carry over these results to a nonlinear situation like the one above, however,
nontrivial modifications are required, which lead to challenging
additional theoretical end experimental investigations.

\bigskip
\section{Summary}
\label{sec:summmary}
In this paper, a TDdG has been developed for a system of Maxwell's equations
with a cubic nonlinearity.
The new capabilities of the proposed method permit that linear and nonlinear
effects of the electric polarization are modeled in an efficient manner
that conserves the energy or is energy stable.
The novel approach allows energy stability both
at the semi-discrete and fully discrete levels, which were not yet available
for the full system of nonlinear Maxwell's equations.
A detailed error estimate is provided for the semi-discrete problem.
The approach is almost completely general and could replace
the electric field formulation, magnetic field formulation, and $A$-formulation.

\end{document}